\documentclass[12pt]{amsart}

\usepackage{enumerate}

\newtheorem{theorem}{Theorem}[section]
\newtheorem{lemma}{Lemma}[section]
\newtheorem{proposition}{Proposition}[section]
\newtheorem{corollary}{Corollary}[section]
\theoremstyle{definition}
\newtheorem{definition}{Definition}[section]
\newtheorem{example}{Example}[section]
\theoremstyle{remark}
\newtheorem{remark}{Remark}[section]
\numberwithin{equation}{section}

\newcommand{\tfa}{time-frequency analysis}

\newcommand{\stft}{short-time Fourier transform}

\newcommand{\tf}{time-frequency}

\newcommand{\tfs}{time-frequency shift}

\newcommand{\modsp}{modulation space}
\newcommand{\psdo}{pseudodifferential operator}

\newcommand{\symbo}{\widetilde M^{\infty, \cA}}

 \def\cF{\mathcal{F}}              
 \def\cS{\mathcal{S}}

 \def\cB{\mathcal{B}}
 
 \def\cG{\mathcal{G}}
 \def\cM{\mathcal{M}}

 \def\cA{\mathcal{A}}

 \def\cL{\mathcal{L}}
 \def\cC{\mathcal{C}}

 \def\cO{\mathcal{O}}
 
 \def\cR{\mathcal{R}}

 \def\cY{\mathcal{Y}}

\def\bR{{\mathbb{R}}}
\def\bZ{{\mathbb{Z}}}
\def\bT{{\mathbb{T}}}
\def\bN{{\mathbb{N}}}

\def\vgf{V_gf}

\def\cca{\Cal C_\ac}
\def\weyl{\sigma^w}

\def\rd{\bR^d}

\def\rdd{{\bR^{2d}}}
\def\zdd{{\bZ^{2d}}}

\def\lrd{L^2(\rd)}
\def\lrdd{L^2(\rdd)}
\def\zd{\bZ^d}
\def\td{\bT^d}

\def\lpm{\ell ^p_v}

\newcommand{\vf}{\varphi}
\def\intrd{\int_{\rd}}
\def\intrdd{\int_{\rdd}}

\def\<{\left<}
\def\>{\right>}

\def\inv{^{-1}}
\def\ud{\, d}

\def\mv1{M_v^1}

\def\Mmpq{M_m^{p,q}}

\def\cF{\mathcal{F}}              
 \def\cS{\mathcal{S}}

 \def\cB{\mathcal{B}}
 
 \def\cG{\mathcal{G}}
 \def\cM{\mathcal{M}}

 \def\cA{\mathcal{A}}

 \def\cC{\mathcal{C}}
 
 \def\cO{\mathcal{O}}
 
 \def\cR{\mathcal{R}}

 \def\cV{\mathcal{V}}

\newcommand{\bba}{{\bf a}}
\newcommand{\bbb}{{\bf b}}
\newcommand{\bbc}{{\bf c}}
\newcommand{\bbd}{{\bf d}}
\newcommand{\bbh}{{\bf h}}
\newcommand{\alb}{\mbox{\boldmath$\alpha$}}
\newcommand{\beb}{\mbox{\boldmath$\beta$}}
\newcommand{\fif}{if and only if}

\newcommand{\bn}{\mathbb{N}}
\newcommand{\br}{\mathbb{R}}

\newcommand{\bz}{\mathbb{Z}}

\newcommand{\Cal}{\mathcal}

\newcommand{\ov}{\overline}
\newcommand{\ch}{\mathbf 1}

\newcommand{\la}{\lambda}
\newcommand{\La}{\Lambda}
\newcommand{\lan}{\langle}
\newcommand{\ran}{\rangle}

\newcommand{\el}{{\ell}^2(\zd)}
\newcommand{\ac}{\Cal A}

\newcommand{\gras}{\cG (\sigma )}
\newcommand{\moda}{\widetilde M^{\infty,\ac}}

\newcommand{\fourd}{\bR ^{4d}}
\newcommand{\phase}{(x,\xi)}
\DeclareMathOperator*{\esup}{ess\,sup\,}

\begin{document}

\title[Banach algebras of pseudodifferential operators]
{Banach algebras of pseudodifferential operators and their almost
diagonalization}

\author{Karlheinz Gr\"ochenig}

\address{ Numerical Harmonic
Analysis Group, Faculty of Mathematics, UNIVERSITY of VIENNA,
Nordbergstrasse 15, A-1090 Wien, Austria}

\email{karlheinz.grochenig@univie.ac.at}

\author{Ziemowit Rzeszotnik}

\address{Mathematical Institute, University of Wroclaw,
Pl.~Grunwaldzki 2/4, 50-384 Wroclaw, Poland}

\email{zioma@math.uni.wroc.pl}

\thanks{The work was supported by the Marie-Curie Excellence Grant
MEXT-CT-2004-517154. It is also part of the  project MOHAWI MA 44
that
is supported by  WWTF (Vienna Science and Technology
Fund)}

\date{\today}
\begin{abstract}
We define new symbol classes for \psdo s and investigate their
pseudodifferential calculus. The symbol classes are parametrized
by
commutative convolution algebras. To every solid
convolution algebra $\cA $ over a lattice  $\Lambda $ we associate
a symbol class
$M^{\infty , \cA } $.
Then  every operator with a symbol in $M^{\infty ,\cA } $ is
almost diagonal with respect to special wave packets (coherent
states
or Gabor frames),  and the rate of almost diagonalization is described
precisely by the underlying convolution algebra $\cA $.
Furthermore,  the corresponding class of \psdo s  is a
Banach algebra of bounded operators on $\lrd $. If a version of
Wiener's lemma holds for $\cA $, then
the algebra of \psdo s is closed under inversion.
The theory contains as a special case the fundamental  results
about
Sj\"ostrand's class and yields a new proof of a theorem of Beals about
the H\"ormander class $S^0_{0,0}$. 
\end{abstract}
\maketitle

\section{Introduction}

We study \psdo s with symbols that are defined by their \tf\
distribution (phase-space distribution). The first symbol class  of
this type  was introduced by Sj\"ostrand~\cite{Sjo94,Sjo95} whose work
has inspired an important line of research by Boulkhemair, Toft, and
others~\cite{Bou97,boul99,Toft01,toft04,toft04a,LM06}. Independently, an alternative
approach with \tf\ methods was developed
in~\cite{GH99,GH03,book,grocomp,gro06}. The starting point of the \tf\
approach  is the observation that the Sj\"ostrand class coincides with
one of the so-called \modsp s that were introduced by Feichtinger
already in 1983~\cite{fei83}. 
The \tf\ approach added several new insights and 
generalizations  to \psdo s with
non-smooth symbols. Recently Sj\"ostrand has again  taken up
the study of \psdo s with non-smooth symbols and substantially
generalized the original definition~\cite{sjo07}. His  motivation is to
study weighted symbol spaces, the boundedness and algebra properties
of the corresponding \psdo s. 

In this  paper we  also study  extensions of the original
Sj\"ostrand class. Our goal is to understand better the following
fundamental questions:
\begin{itemize}
\item Which properties of the
generalized Sj\"ostrand class are responsible for the boundedness of
the corresponding \psdo s on $\lrd $ and on other function spaces?
\item Which properties of the symbol class imply the algebra property of the
operators?
\item Which properties of the symbol class  yield the spectral
  invariance property and thus a
strong form of the functional calculus? 
\end{itemize}

We introduce a family of symbol classes which in general
may contain non-smooth symbols. This family is parametrized by 
 Banach algebras with respect to convolution
on
a lattice $\Lambda \subseteq \rdd$. To each such Banach algebra $\cA $ we
associate a symbol class $\symbo $,  and we  analyze the
properties of the corresponding symbol class. Our main theme is
how properties of the Banach algebra are inherited  by the
corresponding operators, and our results answer the above questions in
the context of the symbol classes $\symbo $. Roughly, the result may
be summarized as follows:

(a) The algebra property of $\cA $ implies that the corresponding 
class of operators is closed under composition. Thus we obtain new
Banach algebras of \psdo s.

(b) If $\cA $ acts boundedly on a solid sequence space $\cY $,
then the
corresponding \psdo s are bounded on a natural function space
associated to $\cY $, a so-called \modsp . In particular, we obtain
the $L^2$-boundedness of \psdo s in this class.  

(c) If $\cA $ is closed under inversion, then the corresponding class
of operators is also closed under inversion. The inverse of a \psdo\ in this
class is again a \psdo\ in this class. This type of result goes back
 to Beals~\cite{beals77} and represents a  strong form of functional calculus.

To be specific, we formulate a special case of  our  main results
explicitly for the 
algebras $\cA = \ell ^\infty _{s}$, which are  defined by the norm
$\|\mathbf{a}\|_{\ell ^\infty _s} = \sup _{k\in \zdd} |\mathbf{a}(k)| (1+|k|)^s
$  (to make $\cA $ into a Banach algebra, we need $s>2d$). In this
case, the corresponding symbol class $\symbo $ is defined by the norm 
\begin{equation}
  \label{eq:cu21}
\|\sigma \|_{\symbo} :=   \sup _{z,\zeta \in \rdd } |(\sigma \cdot
\Phi (\cdot - z) 
  )\,\widehat{}\, (\zeta ) | (1+|\zeta |)^s \, , 
\end{equation}
(where $\Phi$ is the Gaussian) and coincides with the standard \modsp\ $M^{\infty } _{1\otimes
  v_s}$. 

The fundamental result concerns the almost diagonalization of
operators with symbols in $\symbo$ with respect to \tfs s (phase space
shifts).  For $z= (x,\xi )\in \rdd $ let  $\pi (z) f(t) = e^{2\pi i
  \xi \cdot t} f(t-x)$  be the corresponding \tfs . 

\vspace{3 mm}

\textbf{Theorem A.} (Almost diagonalization) \emph{Let $g$ be a nonzero
Schwartz function.  A symbol belongs to the
class $\symbo $ with $\cA = \ell ^\infty _s$, \fif\ } 
$$
|\langle \sigma ^w \pi (z) g, \pi (w) g\rangle | \leq C (1+|w-z|)^{-s}
\qquad \forall w,z\in \rdd \, .
$$

Remarkably, the almost diagonalization of the corresponding \psdo s
characterizes the symbol class completely. 

\vspace{3 mm}

\textbf{Theorem B.} (Boundedness) \emph{If $\sigma \in \symbo$, then
  $\sigma ^w$ is bounded on a whole class of distribution spaces, the
  so-called \modsp s. In particular, $\sigma ^w$ is bounded on  $\lrd
  $.}

\vspace{3 mm}

\textbf{Theorem C.} (Algebra Property) \emph{If $\sigma _1, \sigma _2
  \in \symbo$, then $\sigma _1^w \sigma _2^w = \tau ^w$ for some $\tau
\in \symbo $. Thus the operators with symbols in $\symbo$ for a Banach
algebra with respect to composition. }

\vspace{3 mm}

\textbf{Theorem D.} (Inverse-Closedness) \emph{If $\sigma \in \symbo $
and $\sigma ^w $ is invertible on $\lrd $, then the inverse operator
$(\sigma ^w)\inv = \tau ^w$ possesses again a symbol in $\symbo $. }

\vspace{3 mm}

If the algebra $\cA $ is the convolution algebra $\ell ^1(\zdd )$,
then the corresponding symbol class $\symbo $ coincides with the
Sj\"ostrand class $M^{\infty , 1}$, and we recover the main results of
\cite{Sjo94,Sjo95}. The context of the symbol  classes
$\symbo $ reveals the deeper reasons for why Sj\"ostrand's fundamental
results hold: for instance, the
$L^2$-boundedness of pseudodifferential operators with a symbol in
$M^{\infty ,1}$ can be traced back to the convolution relation $\ell
^1 \ast \ell ^2 \subseteq \ell ^2$. Finally, Wiener's Lemma for
absolutely convergent Fourier series is at the heart of the ``Wiener
algebra property'' of the Sj\"ostrand class. 

The $L^2$-boundedness and the algebra property (Theorems B and C) can
also be derived 
with Sj\"ostrand's methods. In addition, the \tf\ approach yields the
boundedness on a much larger class of function and distribution spaces, the
\modsp s.  These function spaces play the role of smoothness
spaces in \tfa , they may be seen as the analogues of  the Sobolov spaces in the classical
theory. Theorem~D lies deeper and requires 
extended use of Banach algebra concepts. 

As a consequence of these theorems we will give a new treatment of the
H\"ormander class $S^0_{0,0}$. On the one hand, we will characterize
$S^0_{0,0}$ by the almost diagonalization properties with respect to
\tfs s, and on the other hand, we will provide a completely new proof
of Beals' theorem on the inverse-closedness of $S^0_{0,0}$~\cite{beals77}.

As another application of our main results, we will investigate the
action of \psdo s on \tf\ molecules and the almost diagonalization of
\psdo s with respect to \tf\ molecules. This topic is hardly explored
yet, and ours seems to be the first results in \psdo\ theory.

Throughout the paper we use methods from \tfa\ (phase-space analysis)
and Banach algebra methods. We draw on properties of the \stft\ and
the Wigner distribution and the theory of the associated function
spaces, the \modsp s. For the investigation of the almost
diagonalization of \psdo s we will use the well-developed theory of
Gabor frames, which provide a  kind of 
non-orthogonal phase-space expansions of distributions.The main sources are the
books~\cite{folland89,book}.  Another set of tools comes from
the theory of Banach algebras. The inverse-closedness as expressed in
Theorem~D lies quite deep and requires  some new results on Banach
algebras with respect to convolution.    Although the
proof techniques are classical, these arguments seem unusual in the
theory of \psdo s; we therefore  give the proofs in the Appendix.


The paper is organized as follows: In Section~2 we collect the
preliminary definitions of \tfa . Section~3 presents the new  Banach
algebra results that are crucial for the main theorems. In
Section~4 we introduce the generalized Sj\"ostrand classes and show
that the corresponding \psdo s are almost diagonalized with respect to
Gabor frames. In Section~5 we study the boundedness properties of
these operators on \modsp s, the algebra property, and the functional
calculus. Section~6 gives an application to the H\"ormander class
$S^0_{0,0}$ and a new proof of Beals' result in \cite{beals77}. The
final Section~7 is devoted to \tf\ molecules and the almost
diagonalization of operators in the Sj\"ostrand class with respect to
molecules. In the Appendix we give the proofs of the Banach algebra
results of Section~3.

\section{Preliminaries}

In this section we summarize the main definitions and results from
\tfa\ needed here.

\textbf{Time-Frequency Shifts. } We combine time $x\in \rd $ and
frequency $\xi\in \rd $ into a single point $z= (x,\xi ) $ in the
``\tf '' plane $\rdd $. Likewise we combine the operators of
translation and modulation to a \emph{\tfs } and write $$ \pi (z )
f(t) = M_\xi T_x f(t) = e^{2\pi i \xi \cdot t}   f(t-x). $$ The
\emph{short-time Fourier transform (STFT)} of a
function/distribution
  $f$ on $\rd $  with   respect to a window $g$ is defined by
  \begin{eqnarray*}
\vgf (x,\xi) & = & \intrd f(t) \overline{g(t-x)} e^{-2\pi i t\cdot
\xi} \, dt \\ &=&  \langle f, M_\xi T_x g\rangle = \langle f, \pi
(z) g\rangle \, .
  \end{eqnarray*}
The \stft\ of a symbol $\sigma \phase$, for $x,\xi \in \rd$, is a
function on $\fourd $ and will be denoted by $\cV _\Phi \sigma
(z,\zeta )$ for $z,\zeta \in \rdd $ in order to distinguish it
from the STFT of a function on $\rd $.

To compare STFTs with respect to different windows, we will make use of the
pointwise estimate
\begin{equation}
  \label{eq:fi7}
  |V_hf (z)| \leq |\langle k,g\rangle |\inv \, \big( |V_gf|
  \ast |V_hk|\big) (z) \, .
\end{equation}
which holds under various assumptions on $f,g,h,k$, see
e.g.~\cite[Lemma~11.3.3]{book}. 

\textbf{Modulation Spaces.} Let $\varphi (t) = e^{-\pi t\cdot t}$ be
the Gaussian. Then the  modulation space $M^{p,q}_m(\rd)$, $1\leq
p,q\leq \infty $ is 
defined by measuring the norm of the STFT in the weighted space
$L^{p,q}_m(\rdd)$, that is
\[
\|f\|_{M^{p,q}_m(\rd)}=\|V_{\varphi} f\|_{L^{p,q}_m(\rdd)}  =
\biggl( \intrd \biggl(
    \intrd |V_{\varphi} f(x,\xi)|^p\,
m(x,\xi)^p \ud x\biggr)^{q/p} d\xi \biggr)^{1/q} \, .
\]
One of the basic results about \modsp s is the independence of this
definition from the particular test function chosen in the \stft
. Precisely, if the weight satisfies the condition $m(z_1+z_2) \leq C
v(z_1) m(z_2)$ for $z_1,z_2\in \rdd $ and if $g\in M^1_v(\rd )$, i.e,
$\intrdd |V_\varphi g(z)| \, v(z) \, dz <\infty$, then $\|V_gf
\|_{L^{p,q}_m} $ is an equivalent norm on $\Mmpq $~\cite[Thm.~11.3.7]{book}. 

We shall be mainly concerned with the space
$M^{\infty,q}_{1\otimes v} (\rdd)$. For weights $v$ of polynomial
growth ($v(\zeta)=\Cal O(|\zeta|^N)$ for some $N>0$) the space
$M^{\infty,q}_{1\otimes v} (\rdd)$ consists of all $\sigma\in\Cal
S'(\rdd)$ such that the norm
\begin{equation}
  \label{eq:cc21}
  \|\sigma \|_{M^{\infty,q}_{1\otimes v} (\rdd)}
  =\bigg( \intrdd \Big(\esup _{z\in \rdd } |\cV _\Phi
  \sigma (z,\zeta )| \, v(\zeta) \Big)^q  d\zeta \bigg)^\frac 1q
\end{equation}
is finite. For $q=\infty$, the norm is given by
$$  \|\sigma \|_{M^{\infty,\infty}_{1\otimes v} (\rdd)}
  =\esup_{\zeta\in\rdd}\esup _{z\in \rdd } |\cV _\Phi
  \sigma (z,\zeta )| \, v(\zeta).
$$
For the specific weight
$v_s(\zeta)=\lan\zeta\ran^s=(1+|\zeta|^2)^{\frac s2}$, $s\ge 0$, the
space  $M^{\infty,q}_{1\otimes v_s} (\rdd)$ shall be denoted by
$M^{\infty,q}_{1\otimes \lan\cdot\ran^s} (\rdd)$. Also, the space $M^{p,p}_m$ shall be denoted by $M^{p}_m$.


\textbf{Weyl Calculus.} The \emph{Wigner distribution} of $f,g \in
\lrd $ is defined as
\begin{equation*}
  \label{eq:1a}
    W(f,g)(x,\xi)=\int_{\rd} f\Big(x+\frac{t}2\Big)
    \overline{g\Big(x-\frac{t}2\Big)} e^{-2\pi
    it\cdot\xi }\,dt.
\end{equation*}

 The  Weyl transfrom $\sigma^w $ of a symbol $\sigma\in \Cal S ' (\rdd
 )$ is defined by the sesquilinear form
 \begin{equation}
   \label{eq:f2}
   \langle \sigma ^w f, g \rangle =  \langle \sigma , W(g,f)\rangle
   \qquad f, g \in \cS (\rd ) \, .
 \end{equation}
Usually, the Weyl correspondence is defined by the integral
operator
\begin{equation}
  \label{eq:f2a}
\sigma ^w f(x) = \intrd \sigma \Big(\frac{x+y}{2},\xi \Big)
e^{2\pi i (x-y)\cdot \xi } f(y) \, dy d\xi \, ,
\end{equation}
but this definition is somewhat restrictive, and we will not need
this
particular formula.

The composition of two Weyl transforms defines a twisted product
between symbols via
$$
\sigma ^w \tau ^w = (\sigma \, \sharp \, \tau )^w \, ,
$$
Usually the analysis of the twisted product is based on the formula~\cite{folland89,HorIII85}
\begin{equation}
  \label{eq:ll1}
(  \sigma \, \sharp \tau ) (x,\xi ) = \intrd \intrd \intrd \intrd
\sigma (u,\zeta ) \tau (v,\eta ) \, e^{4\pi i [(x-u)\cdot (\xi -\eta )
  - (x-v)\cdot (\xi - \zeta )]} \, dudvd\eta d\zeta \, ,
\end{equation}
but  it is a distinctive feature of our approach that we will not
need the explicit formula~\eqref{eq:ll1}.

\section{Discrete Banach Algebras}

In this section we present the necessary Banach algebra methods.
We
have not found them in the literature, and will give the proofs in
the
appendix.

Throughout this section $\ac$  denotes a solid involutive Banach
algebra with respect to convolution and indexed by a discrete
subgroup $\Lambda $ of $\rd $ (in the remaining parts of the paper
$\La\subset\rdd$). The elements of $\cA $ are sequences
$\bba (\lambda ), \lambda \in \Lambda $, where $\Lambda = A\zd $ is
a discrete subgroup of full rank (a lattice) in $\rd $ (thus $\det
A \neq 0$). The involution is defined as $\bba^*(\lambda ) =
\overline{\bba(-\lambda )}$. The norm of $\cA $ satisfies the
usual inequalities $\|\bba ^* \|_{\ac } \leq \|\bba \|_{\ac }$ and
\begin{equation}\label {norm}
\|\bba \ast \bbb \|_\ac\le\|\bba \|_\ac\|\bbb\|_\ac \qquad \text{
for
  all } \, \bba , \bbb \in \ac \, .
\end{equation}
Furthermore, the solidity of $\cA $ says that if $|\bba (\lambda )
|\le |\bbb (\lambda ) |$ for all $\lambda \in \Lambda $  and
$\bbb \in\ac$, then also  $\bba \in\ac$ and $\|\bba\|_\ac\le\|\bbb
\|_\ac$.

Due to the solidity, we may assume without loss of generality that
the standard sequences $\delta _\lambda $ defined by $\delta
_\lambda (\mu ) = 1$ for $\lambda = \mu $ and $\delta _\lambda
(\mu ) = 0$ for $\lambda \neq \mu $ belong to $\cA $, otherwise we
switch to the sublattice $\Lambda ' = \{ \lambda \in \Lambda :
\delta _\lambda \in \cA \}$.

The solidity is a strong assumption, as is demonstrated by the
following result.

\begin{theorem}\label{alla}
Let $\Cal A$ be an involutive  Banach algebra with respect to
convolution  over  a lattice $\La\subset\rd$. If $\ac$ is solid,
then $\Cal A $ is
continuously embedded in $
\ell ^1(\La)$.
\end{theorem}

Thus $\ell ^1(\Lambda )$ is the maximal solid convolution algebra
on
the discrete group $\Lambda $.

\begin{example}\label{example1} Let $v$ be a non-negative function on $\Lambda
$
and let $\ell^q_v (\Lambda )$ be defined by the norm
$\|\mathbf{a}\|_{\ell^q_v } =  \| \mathbf{a} v\|_{q}$. Then $\ell
^1_v(\Lambda )$ is a Banach algebra \fif\ $v(\lambda +\mu ) \leq
C\, v(\lambda ) v(\mu ), \lambda ,\mu \in \Lambda $ ($v$ is
submultiplicative), and $\ell ^\infty _v(\Lambda )$ is a Banach
algebra, \fif\ $v\inv \ast v\inv \leq C v\inv $ ($v$ is
subconvolutive)~\cite{feichtinger79}. The standard family of weights is given
as $\lan\la\ran^s$, where $\lan\la\ran=\sqrt{1+|\la|^2}$. The
corresponding solid Banach algebras are $\ell
^1_{\lan\cdot\ran^s}(\Lambda )$, for $s\ge 0$ and $\ell
^\infty_{\lan\cdot\ran^s}(\Lambda )$, for $s> d$.

Another example of a solid Banach algebra is the so-called Krein
algebra $\ell ^1(\La) \cap \ell ^2 _{d/2}(\La)$ with the norm
$\|\mathbf{a}\|_1 + \Big( \sum _{\lambda \in \Lambda }  |a_\lambda
|^2 |\lambda |^{d} \Big) ^{1/2}$, see \cite{bra75}.

 For certain weight
sequences  $v$ the weighted $\ell ^q$-space $\ell^q_v(\La) $ is
also a convolution algebra, see~\cite{feichtinger79}.

\end{example}

 To every  solid Banach
algebra under convolution we can attach a Banach algebra of
matrices.
Indeed, using an idea of Baskakov~\cite{Bas90}, we  define an algebra
of matrices that are ``dominated by convolution operators in $\cA
$''.

\begin{definition}\label{bas}
Let  $A $ be a matrix on $\Lambda $  with entries $a_{\lambda
\mu}$, for $\lambda , \mu \in \Lambda $, and  let $\mathbf{d}_A$
be the sequence with entries $\bbd_A(\mu)$ defined by
 \begin{equation}
\label{eq:ma1} \bbd_A (\mu ) = \sup _{\lambda \in \Lambda } |a
_{\lambda , \lambda  -\mu}| \, .
\end{equation}
We say that the matrix $A $ belongs to $\cca $, if  $\bbd_A$
belongs to $\cA $. The norm in $\cca $  is given by
\begin{equation}
\label{eq:ma2} \|A\|_{\cca } = \|\bbd_A\|_{\cA } \, .
\end{equation}
\end{definition}

Note that $\bbd_A(\mu)$ is the supremum of the entries in the
$\mu$-th diagonal of $A$, thus the $\cca $-norm describes a form
of the off-diagonal decay of a matrix.

We first list some  elementary properties of $\cca $.

\begin{lemma}\label{lem:ca1}
Assume that $\cA $ is a solid Banach algebra under
convolution. Then:

(i) $\cca $ is a Banach algebra under matrix multiplication (or
equivalently, the composition of the associated operators).

(ii) Let $\cY$ be a solid Banach space of sequences on $\Lambda $.
If
$\cA $ acts boundedly on a solid space $\cY$ by
convolution ($\cA \ast \cY \subseteq \cY$), then $\cca $ acts
boundedly on
$\cY$, i.e.,
\begin{equation}
\label{eq:ma3} \|A\mathbf{c}\|_\cY \leq  \|A\|_{\cca }
\|\mathbf{c}\|_{\cY} \qquad \text{ for all } \, A\in \cca,
\mathbf{c}\in \cY \, .
\end{equation}

(iii) In particular, since $\cA \subseteq \ell ^1(\La)$, we may
identify $\cca  $ as a (Banach) subalgebra of $\cB (\ell
^2(\La))$.
\end{lemma}

\begin{proof}
(i) Let $A= (a_{\la\mu}), B= (b_{\la\mu})\in \cca $, then by
definition $|a_{\la\mu}| = |a_{\la,\la-(\la-\mu)}| \leq d_A
(\la-\mu) $. Consequently matrix multiplication is dominated by
convolution in the sense that
\begin{eqnarray*}
|(AB)_{\la, \la-\mu}| &\leq &\sum _{\nu\in \La } |a_{\la\nu} | \,
|b_{\nu,\la-\mu}|
\\&\leq &  \sum _{\nu\in \La } \bbd_A( \la-\nu)  \, \bbd_B(\nu-(\la-\mu))      \\
&=& (\bbd_A \ast \bbd_B)(\mu) \, ,
\end{eqnarray*}
and so $\bbd_{AB}(\mu) = \sup _{\la\in \La } |(AB)_{\la,\la-\mu}|
\leq (\bbd_A \ast \bbd_B)(\mu)$. Since $\bbd_A, \bbd_B\in \cA $
and $\cA $ is a Banach algebra under convolution, we find that

$$
 \|AB\|_{\cca }  = \|\bbd_{AB}\|_{\cA } \leq \|\bbd_A \ast
\bbd_B\|_{\cA } \leq \|\bbd_A \|_{\cA } \, \| \bbd_B \|_{\cA } =
\|A\|_{\cca }\, \|B\|_{\cca } \,. $$ (ii) Since $|a_{\la\mu}| \le
\bbd_A (\la-\mu) $, we obtain the pointwise inequality

$$ |A\mathbf{c}(\la)|   = \Big|\sum _{\mu\in \La } a_{\la\mu}
\bbc(\mu) \Big| \leq \sum _{\mu\in \La } \bbd_A (\la-\mu)
|\bbc(\mu)| = \big( \bbd_A \ast |\mathbf{c}|\big) (\la) \, . $$

Using the hypothesis on $\cY$, we conclude that
\begin{equation*}
\|A\mathbf{c}\|_\cY \leq \|\bbd_A \ast |\mathbf{c}| \, \|_\cY \leq
\|\bbd_A \|_{\cA } \, \| \, |\mathbf{c}| \, \|_\cY  = \|A\|_{\cca
} \, \|\mathbf{c} \, \|_\cY  \, ,
\end{equation*}
since $\cY$ is solid and $\cA $ acts on $\cY$ by convolution.

(iii) follows by choosing $\cY=\ell ^2(\La )$ and by Young's
inequality $\ell ^1(\La) \ast \ell ^2(\La) \subseteq \ell
^2(\La)$.
\end{proof}

Whereas $\cA $ is a commutative Banach algebra,  $\cca $ is highly
non-commutative, the transition from $\cA $ to $\cca $  can be
thought of as a non-commutative extension of convolution algebras
of sequences on $\La $.

One of the main questions about the matrix algebra $\cca $ is
whether
the inverse of a matrix in $\cca $ is again in $\cca $, or in
other
words, we ask whether the off-diagonal decay  described by $\cA $
is
preserved under inversion.

We recall the following definition. Let $\cA \subseteq \cB $ be
two Banach algebras with a common unit element. Then $\cA $ is
\emph{inverse-closed} in $\cB $, if $\bba\in \cA $ and $\bba\inv
\in \cB
$
implies that $\bba\inv \in \cA $.

In the following, we identify an  element $\bba\in \cA \subseteq
\ell ^1(\Lambda )$ with the corresponding convolution operator
$C_\mathbf{a} \mathbf{b} = \mathbf{a} \ast \mathbf{b}$. In this
way, we may treat $\cA $ as a Banach subalgebra of $\cB (\ell
^2(\Lambda ))$, the algebra of  bounded operators on $\ell
^2(\Lambda )$.

The following theorem gives a complete characterization  of  when
the non-commutative extension $\cca $ is inverse-closed in
 $\cB (\ell ^2(\La ))$.  For $\cA = \ell ^1_v(\zd)$ and
 $\cA = \ell ^\infty _v(\zd)$ this
characterization is due to 
 Baskakov~\cite{Bas90}. 
Our formulation is new and reveals more clearly what the main
conditions are. 
Recall that the spectrum $\widehat{\cA }$ of a
commutative Banach algebra $\cA $ consists of all multiplicative
linear functionals on $\cA $. We denote the standard basis of $\ell
^2(\Lambda )$ by $\delta _\lambda , \lambda \in \Lambda $. 

The following theorem is crucial for the functional calculus of
\psdo s.

\begin{theorem} \label{bask}
Assume that $\cA $ is a solid convolution algebra of sequences on
a lattice $\Lambda\subset \rd $. Then the following are
equivalent:
\begin{itemize}

  \item[(i)] $\cA $ is inverse-closed in $\cB (\ell^2)$.

  \item[(ii)] $\cca $ is inverse-closed in $\cB (\ell^2)$.

  \item[(iii)] The spectrum $\widehat{\cA }\simeq \bT^d $.

  \item[(iv)] The weight $\omega (\lambda ) = \|\delta _\lambda \|_\cA
    $ satisfies the GRS-condition $\lim _{n\to \infty } \omega
    (n\lambda )^{1/n} = 1$ for all $\lambda \in \Lambda $. 
  \end{itemize}
\end{theorem}
The main point of condition (iv) is that we have an easy condition to
check whether the non-commutative matrix algebra $\cC _{\cA }$ is
inverse-closed in $\cB (\ell ^2)$. 
We defer the proof of this theorem to the appendix, because  this
paper is on \psdo s and Banach algebras are only a tool. The proof
is
an extension and re-interpretation of Baskakov's argument.

\begin{example}\label{example2} By a theorem of Gelfand, Raikov, and
Shilov~\cite{gelfandraikov}, the weighted  convolution algebra $\ell ^1_v (\La
)$ is inverse-closed in $\cB (\ell ^2(\La))$, \fif\ the weight $v$
satisfies the (so-called GRS) condition $\lim _{n\to \infty }
v(n\lambda )^{1/n} =1$ for all $\lambda \in \Lambda $. Thus
Theorem~\ref{bask} implies that the non-commutative matrix algebra
$\cC _{\ell ^1_v(\La)}$ is inverse-closed in $\cB (\ell ^2(\La))$,
\fif\ the weight satisfies the GRS-condition. This is the result of
Baskakov~\cite{Bas90}. Since the polynomial  weight
$\lan\la\ran^s$ satisfies the GRS-condition when $s\ge 0$, the
algebras $\ell^1_{\lan\cdot\ran^s}(\La)$ and
$\cC_{\ell^1_{\lan\cdot\ran^s}(\La)}$ are inverse-closed for $s\ge
0$.
Similarly,  the algebras $\ell^\infty_{v}(\La)$ and
$\cC_{\ell^\infty_{v}(\La)}$ are inverse-closed if $v$ satisfies
the GRS-condition and $v^{-1}\in \ell^1(\La)$. Therefore the
standard algebras $\ell^\infty_{\lan\cdot\ran^s}(\La)$ and
$\cC_{\ell^\infty_{\lan\cdot\ran^s}(\La)}$ (the Jaffard class) are
inverse-closed for $s>d$, see~\cite{jaffard90,Bas90,GL04}.
\end{example}

As a consequence of Theorem~\ref{bask} we draw the following.

\begin{corollary}[Spectral Invariance]\label{specin}
Assume that $\widehat{\cA } \simeq \bT ^d$. Then
\begin{equation}\label{eq:ca11}
 \mathrm{Sp}_{\cB (\ell ^2)} (A) = \mathrm{Sp}_{\cca } (A) \qquad \text{
for all } \, A \in \cca \, .
\end{equation}
Moreover, if $\cA $ acts boundedly  on a solid sequence space
$\cY$,
then
\begin{equation}\label{eq:ca12}
\mathrm{Sp} _{\cB (\cY )} (A) \subseteq \mathrm{Sp}_{\cB (\ell ^2)} (A) =
\mathrm{Sp}_{\cca } (A) \qquad \text{ for all } \, A \in \cca \, .
\end{equation}
\end{corollary}

\begin{proof}
The spectral identity \eqref{eq:ca11} is just a reformulation of the
fact that $\cC _{\cA } $ is inverse-closed in $\cB (\ell ^2)$. If
$A\in \cC _{\cA }$ is invertible on $\ell ^2$, then also $A\inv \in
\cC _{\cA }$. Thus by Lemma~\ref{lem:ca1} $A\inv $ is bounded on
$\cY$. Consequently, if $\lambda \not  \in \mathrm{Sp}_{\cB (\ell
  ^2)}(A)$, then $\lambda \not  \in \mathrm{Sp}_{\cB (\cY )}$, which
is the inclusion~\eqref{eq:ca12}.
\end{proof}

\section{Almost diagonalization}

In this section we introduce the general symbol classes $\symbo $ that
are parametrized by a solid convolution algebra $\cA $. The main
result explains how \psdo s  with symbols in $\symbo $ are almost
diagonalized. 

In order to define the symbol classes we recall
the basic information about amalgam spaces and Gabor frames. The
definitions are adapted to our needs, for the general theory of
amalgam spaces we refer to ~\cite{fournier-stewart85,feichtinger90}, for Gabor frames
to~\cite{book,chr03}. 

\textbf{Amalgam spaces.}  Fix  a  Banach algebra $\ac$ of
sequences on a lattice $\Lambda\subseteq \rdd  $ and a  relatively compact
fundamental domain $C$ containing the origin. We     say that a
function $F\in L^\infty_{\text {loc}}(\rdd)$ belongs to  the associated
amalgam space $W(\ac)$, 
if the sequence $\mathbf{a}$ of local suprema
\[
\bba(\la)=\esup _{\zeta\in\la+C} F(\zeta)
\]
belongs to $\ac$.  The norm on $W(\ac)$ is given by
$\|F\|_{W(\ac)}=\|\mathbf{a}\|_\ac$. For completeness we mention that
this definition is independent of the lattice $\Lambda $ and the fundamental
domain $C$. We will use frequently that   $W(\cA )$ is an involutive 
Banach algebra with respect to convolution on $\rdd
$~\cite{feichtinger83}.  

\textbf{Gabor Frames.} Let $\La=A\zdd$ be a lattice with $|\det
A|<1$. We say that a set $\Cal
G(g,\La)=\{\pi(\la)g:\la\in\La\}$ forms a tight Gabor frame (with constant 1) for
$\lrd$ if 
\begin{equation}
  \label{eq:ll2}
  \|f\|_2^2 = \sum _{\lambda \in \Lambda } |\langle f, \pi (\lambda
  ) g\rangle |^2 \qquad  \text{ for all } f\in\lrd \, .
\end{equation}
As a consequence every $f\in \lrd $ possesses the  {\it tight
frame expansion}
\[
f=\sum_{\la\in\La}\lan f,\pi(\la)g\ran\pi(\la)g\, ,
\]
with unconditional convergence in $\lrd $.

In our considerations it will be important to use  a frame
with a window $g$ that satisfies an additional assumption $V_gg \in
W(\ac)$. We will take the existence of a tight frame with $V_gg \in W(\ac
)$ for granted and will not worry about the subtle existence problem. It is
known that  tight Gabor frames with $g\in \lrd $ exist  for every
lattice $\Lambda = A\zdd $ with $|\det A| <1$~\cite{bekka04}.   
If $\Lambda = \alpha \zd \times \beta \zd $ and $\alpha \beta <1$ (or
the symplectic image of such a lattice), then for every $\delta>0, 0<b<1$
there exist tight Gabor 
frames with $|V_gg (z) | \leq  Ce^{-\delta|z|^b}$(e.g.,~\cite{GL04}).


\textbf{Symbol Classes.}  Now we define a new family of symbol classes. As in
  ~\cite{grocomp}  the  grand symbol
$\gras$ of a symbol $\sigma $  is given by
\[
\gras(\zeta):=\esup _{z\in \rdd } |\cV _\Phi \sigma (z,\zeta )|\, ,
\]
where $\Phi(z) = e^{-\pi z\cdot z/2} $ is the Gaussian or some equivalent window (see below). Furthermore, set 
$j(\zeta)=(\zeta_2,-\zeta_1)$ for $\zeta=(\zeta_1,\zeta_2)\in  \rdd $.

\begin{definition}\label{def:sjocl}  A symbol $\sigma$ belongs to the generalized
  Sj\"ostrand class 
$\widetilde M^{\infty,\ac}(\rdd )$ if $\gras\circ j\in W(\ac)$.
The norm on $\widetilde M^{\infty,\ac}$ is given
by 
$$\|\sigma\|_{\widetilde M^{\infty,\ac}}= \|\gras\circ j\|_{
W(\ac)}\, .
$$
\end{definition}

 The tilde over $M$ indicates that we take into
account the mapping $j$. (It could be avoided by using the symplectic
Fourier transform in addition to the ordinary Fourier transform.)

\begin{remark}
  If $\cA = \ell^q_v(\Lambda)$,
where the weight $v$ is continuous and submultiplicative on
$\rdd$, then   we obtain the  standard \modsp s 
\begin{equation}\label{lpm}
\widetilde M^{\infty,\ac}= M^{\infty,q}_{1\otimes v\circ
j^{-1}}(\rdd), \quad 1\le q\le \infty
\end{equation}
In particular, if $\cA = \ell ^1(\Lambda )$, then $\symbo $ is the
Sj\"ostrand class $M^{\infty ,1}$. 
Thus the symbol class $\widetilde M^{\infty,\ac}$ is a generalized
modulation space. To prove (\ref{lpm}) we observe that, by
\eqref{eq:cc21}, $\sigma\in M^{\infty,q}_{1\otimes v\circ
j^{-1}}(\rdd)$,  \fif\ $\gras\circ j\in L^q_{ v}(\rdd)$. Since $v$ is
continuous and submultiplicative, the inclusion 
 $W(\ell^q_v(\Lambda))\subseteq L^q_{ v}(\rd)$ holds for every
 lattice $\La\subset \rdd$. This shows that
$\widetilde M^{\infty,\ac}\subseteq M^{\infty,q}_{1\otimes v\circ
j^{-1}}(\rdd)$. Conversely,  if
$\sigma\in M^{\infty,q}_{1\otimes v\circ j^{-1}}(\rdd)$, then
$\gras\in W(\ell^q_{v\circ j^{-1}}(\La))$ for every 
lattice $\La\subset \rdd$ by  \cite[Thm~12.2.1]{book}, and we have
equality in~\eqref{lpm}.
\end{remark}

To consolidate Definition~\ref{def:sjocl}, we need to establish its 
independence of the particular window $\Phi $. This is done as in the
case of \modsp s in~\cite[Ch.~11.3,11.4]{book}, but requires a few adjustments. 

\begin{lemma}\label{lemind}
Let $\Psi(z) = e^{-\pi z\cdot z/2}$ be the Gaussian. If a window $\Phi \in \lrdd $ satisfies the condition
\begin{equation}
  \label{eq:fi9}
  F(\zeta ) := \intrdd |\cV _\Psi \Phi (z,j(\zeta ))| dz \in W(\cA )
  \, ,
\end{equation}
then the definition of $\symbo $ does  not depend on the window.   
\end{lemma}

\begin{proof}
  We use the pointwise estimate~\eqref{eq:fi7} in the form
$$
|\cV _\Psi \sigma (z,\zeta )| \leq \frac{1}{\|\Phi \|_2^2}  \intrdd |\cV _\Phi
\sigma (z-u, \eta  ) | \, |\cV _\Psi \Phi (u, \zeta -\eta) | \,du
d\eta \, .
$$
After taking the supremum over $z$ and inserting \eqref{eq:fi9},  we
obtain
$$
\Cal G_\Psi (\sigma ) (\zeta ) \leq \Big( \Cal G_\Phi (\sigma )  \ast (F\circ j\inv)
\Big)(\zeta ) \, ,
$$
where $\Cal G_\Psi $ is the grand symbol with respect to the window $\Psi
$. 
Consequently, since $W(\cA )$ is a Banach algebra, we obtain
$$
\|\Cal G_\Psi (\sigma )\circ j \|_{W(\cA )} \leq \|\Cal G_\Phi (\sigma )\circ j\|_{W(\cA ) } \,
\|F\|_{W(\cA    )}\, .$$
By interchanging the roles of $\Phi $ and $\Psi $, we obtain the norm
equivalence
$$
\|\Cal G_\Phi (\sigma )\circ j \|_{W(\cA )} \asymp \|\Cal G_\Psi (\sigma ) \circ j\|_{W(\cA )}
= \|\sigma \|_{\symbo } \, .
$$
\end{proof}

Next, to establish the link between pseudodifferential operators and
the generalized Sj\"ostrand classes $\symbo $, we will need windows of
the form $\Phi = W(g,g)$ for a suitable function $g$ on $\rd
$. According to Lemma~\ref{lemind} we need to determine a class of
functions $g$ such that $\Phi = W(g,g)$ satisfies
condition~\eqref{eq:fi9}. This is explained in the following lemma.

\begin{lemma}
  \label{lem:htht}
If $V_gg \in W(\cA )$, then $\Phi = W(g,g)$ satisfies condition
\eqref{eq:fi9}. 
\end{lemma}

\begin{proof}
We first show that the hypothesis implies that $V_\vf g \in W(\cA )$. 
Assume first that $\langle g, \vf \rangle \neq 0$, then by
\eqref{eq:fi7}
$$
|V_\vf g| \leq  |\langle g,\vf \rangle |\inv \, \big(|V_gg| \ast
|V_\vf \vf|\big) \, ,
$$
thus $\|V_\vf g\|_{W(\cA )} \leq |\langle g,\vf \rangle |\inv \,
\|V_gg\|_{W(\cA )} \, \|V_\vf \vf \|_{W(\cA )}<\infty $. If $\langle g, \vf
\rangle = 0$, then choose a $\vf _1$ (e.g., a linear combination of
Hermite 
functions), such that $\langle \vf _1, g\rangle \neq 0$ and $\langle \vf _1,
\vf \rangle \neq 0$ and apply the above argument twice.

Next  observe that $e^{-\pi z\cdot z/2} = W(\vf , \vf )(z)$,
$z\in \rdd $, where $\vf (t) = e^{-\pi t\cdot t}$,  and we use the
``magic formula'' for the STFT of a 
Wigner distribution~\cite[Lemma~14.5.1]{book}:
$$
|\cV _{W(\vf , \vf )} W(g,g)(z,\zeta )| =|V_\vf g (z+\frac{j(\zeta )
  }{2})  |\, |V_\vf g (z-\frac{j(\zeta )}{2})    | \, .
$$
Consequently,
$$
 F(\zeta ) = \intrdd  |\cV _{W(\vf , \vf )} W(g,g)(z,j(\zeta ))| \, dz = \big(|V_\vf g|
\ast |(V_\vf g)^*|\big) (\zeta )
$$
where $|(V_\vf g )^* (z)| =|V_\vf g (-z)|=|V_g \vf (z)|$. 
Since both  $V_\vf g \in W(\cA )$ and   $V_g \vf  \in W(\cA )$, we
have $F \in W(\cA )$, as claimed. 
\end{proof}

\textbf{Almost diagonalization.} 
An important
identity that establishes a relation between the \stft\ of a
symbol $\sigma $ and its  Weyl transform $\sigma ^w$ is the following.

\begin{lemma}[Lemma~3.1\cite{gro06}]
  Assume that
  $\sigma \in M^\infty (\rdd )$ and $g\in M^1 (\rd )$. 
  If we choose the window
  $\Phi $ to be the Wigner distribution $\Phi  =
  W_gg$, then
\begin{equation}\label{eq:4a}
\big|\langle \sigma ^w \pi (z)g , \pi (w) g \rangle \big| =
\Big|\cV_\Phi \sigma \Big(\frac{w+z}{2}, j(w-z)\Big)\Big|
\end{equation}
 for all $w,z\in \rdd $.
\end{lemma}

 Read backwards, this  formula yields
\begin{equation}
  \label{eq:5a}
  |\cV_\Phi \sigma (u,v)| = \Big|\Big\langle \sigma ^w  \pi \Big(u-\frac{j\inv
  (v)}{2}\Big)g, \pi \Big(u+\frac{j\inv   (v)}{2}\Big)g \Big\rangle
\Big| \, .
\end{equation}

Both formulas \eqref{eq:4a} and \eqref{eq:5a} hold pointwise. 

In order
to check the membership of a symbol to $\symbo  $,
we need at least that the window $\Phi = W(g,g)$ used to measure the
$\symbo$-norm  satisfies
condition~\eqref{eq:fi9}. This is guaranteed by
Lemma~\ref{lem:htht}. In the remainder of the paper the condition
$V_gg \in W(\cA )$ will thus be the standing assumption of $g$. In
particular, we may use formulas \eqref{eq:4a} and \eqref{eq:5a},
because $V_gg \in W(\cA )$ implies that $g\in M^1$ and $\symbo
\subseteq M^\infty $. 

 Our main theorem is a far-reaching extension of
Theorem 3.2 in~\cite{gro06}.

\begin{theorem}\label{aldia}
Let $\ac$ be a solid Banach algebra with involution on a lattice
$\Lambda \subseteq \rdd $. Fix  a window $g$ such that $V_g g\in
W(\ac)$ and let  $\Cal G(g,\Lambda)$ be a tight Gabor frame for
$L^2(\rd)$. Then the following are equivalent for a distribution
$\sigma \in \cS ' (\rdd )$.
\begin{itemize}
\item[(i)] $\sigma\in \widetilde M^{\infty,\Cal A}$.
\item[(ii)] There exists a function $H\in W(\ac)$
such that
\begin{equation} \label{eq:6}
|\lan \sigma^{w}\pi(z)g,\pi(w)g\ran|\le H(w-z) \quad \text{ for
}\, w,z\in\rdd.
\end{equation}
\item[(iii)] There exists a sequence $\bbh\in\Cal A$ such
that
\begin{equation} \label{eq:6b}
|\lan \sigma^{w}\pi(\mu)g,\pi(\lambda)g\ran|\le \bbh(\lambda-\mu)
\quad \text{ for }\, \lambda,\mu\in\Lambda.
\end{equation}
\end{itemize}
\end{theorem}
\begin{proof}
 The equivalence of (i) and (ii) follows easily by using formulae
 (\ref{eq:4a}) and (\ref{eq:5a}).

(i) $ \Rightarrow  $ (ii). By (\ref{eq:4a}) we have
\begin{equation}\label{eq:7}
\big|\langle \sigma ^w \pi (z)g , \pi (w) g \rangle \big| =
\Big|\cV_\Phi \sigma \Big(\frac{w+z}{2}, j(w-z)\Big)\Big|\le
\gras(j(w-z)).
\end{equation}
Since $\gras\circ j\in W(\ac)$, we may take $H=\gras\circ j$ as
the dominating function in (\ref{eq:6}).

(ii) $ \Rightarrow  $ (i). By (\ref{eq:5a}) and (\ref{eq:6}) we
have
\begin{equation}\label{eq:8}
\gras(j(\zeta))=\esup_{z\in\rdd}|\cV_\Phi \sigma (z,
j(\zeta))|=\esup_{z\in\rdd}\Big|\Big\langle \sigma ^w  \pi
\Big(z-\frac{\zeta}{2}\Big)g, \pi \Big(z+\frac{\zeta}{2}\Big)g
\Big\rangle \Big| \le H(\zeta).
\end{equation}
Thus, $\gras\circ j\le H$. Since $H\in W(\ac)$ and $\ac$ is solid
we get that $\gras\circ j\in W(\ac)$.

(i) $ \Rightarrow  $ (iii). This implication follows 
from formulae (\ref{eq:4a}) as well. Indeed, for
$\la,\mu\in\Lambda$ we have, as in (\ref{eq:7}), that
\[
\big|\langle \sigma ^w \pi (\mu)g , \pi (\la) g \rangle \big| \le
\gras(j(\la-\mu)).
\]
Let $\bbh(\la)=\esup _{\zeta\in\la+C} \gras(j(\zeta))$. Since the
fundamental domain  $C$
is assumed to contain the origin,  we have  $\gras(j(\la))\le
\bbh(\la)$. Therefore, the above inequality can be extended to
\begin{equation}\label{eq:8a}
\big|\langle \sigma ^w \pi (\mu)g , \pi (\la) g \rangle \big| \le
\gras(j(\la-\mu))\le \bbh(\la-\mu).
\end{equation}
Since $\sigma\in \widetilde M^{\infty,\Cal A}$ means that
$\bbh\in\ac$, we obtain (\ref{eq:6b}).

(iii) $ \Rightarrow  $ (ii). This  implication is more technical.
Only
here we use  the
assumption that $\Cal G(g,\Lambda)$ is a tight Gabor frame for
$L^2(\rdd )$. Consider  the tight frame expansion of $\pi(u)g$
\begin{equation}
  \label{eq:12}
\pi (u)g = \sum _{\nu \in \Lambda } \langle \pi (u)g,\pi (\nu ) g
\rangle \pi (\nu ) g \, ,
\end{equation}
for every $u\in\rdd  $. Since we assume that $V_g g\in W(\ac)$,
 the sequence $\alb$ of local suprema
\begin{equation}
  \label{eq:13b}
\alb (\nu ) = \sup _{\zeta\in C} |V_g g(\nu +\zeta)| = \sup
_{\zeta\in C} |\langle g,\pi (\nu+\zeta) g \rangle |=\sup
_{\zeta\in C} |\langle \pi (\zeta)g,\pi (-\nu) g \rangle | \, ,
\quad \nu \in \Lambda \, ,
\end{equation}
belongs to $\ac $.

For given $z,w \in \rdd $,  we write them uniquely as $z=\mu +u'$,
$w=\lambda + u$, where  $\lambda ,\mu \in \Lambda$ and $u,u' \in
C$. Inserting  the expansions~\eqref{eq:12} and the definition of
$\alb$ in the matrix entries,  we find that (\ref{eq:6b}) yields
the following estimate
\begin{eqnarray*}
\lefteqn{  |\langle \sigma ^w \pi (\mu +u')g, \pi (\lambda
+u)g\rangle | =
  |\langle \sigma ^w \pi (\mu ) \pi (u')g, \pi (\lambda) \pi
  (u)g\rangle |} \\
&\leq & \sum _{\nu, \nu' \in \Lambda } |\langle \sigma ^w \pi (\mu
+\nu ')g, \pi (\lambda +\nu )g\rangle| \, |\langle \pi (u')g, \pi
(\nu') g \rangle | \, |\langle \pi (u) g, \pi (\nu )g \rangle | \\
&\leq & \sum _{\nu, \nu '\in \Lambda } \bbh(\lambda- \mu +\nu
-\nu ') \alb (-\nu') \alb (-\nu ) \\ &=& (\bbh\ast \alb \ast
\alb^* )(\lambda -\mu )\, ,
\end{eqnarray*}
where  $\alb ^*(\lambda ) = \overline{\alb (-\lambda )}$ is the
involution on $\ac$. Since $\alb\in\ac$, $\alb^*\in\ac$ and
$\bbh\in\ac$, we see that $\bbh\ast \alb \ast \alb^*\in \ac$ as
well.

To find the dominating function $H$ postulated in (\ref{eq:6}),
we set
\begin{equation}\label{eq:14}
H(\zeta) = \sum _{\nu \in \Lambda } (\bbh\ast \alb  \ast \alb ^* )
(\nu ) \,  \chi _{C-C}(\zeta-\nu) \, ,
\end{equation}
where $\zeta\in\rdd$. 
Our previous estimate says that 
\begin{eqnarray*}
\lefteqn{|\langle \sigma ^w \pi (z)g, \pi (w)g\rangle |\le
(\bbh\ast \alb \ast \alb^*   )(\lambda -\mu )=(\bbh\ast \alb \ast
\alb^* )(\lambda -\mu )\chi _{C-C}(u-u')}\\ &\le &\sum _{\nu \in
\Lambda } (\bbh\ast \alb  \ast \alb ^* ) (\nu ) \,  \chi
_{C-C}(\la-\mu+u-u'-\nu)=H(\la-\mu+u-u')=H(w-z)\, ,
\end{eqnarray*}
and $H$ satisfies \eqref{eq:6}. 

To finish, we need to show that $H\in W(\ac)$. Let
$\beb(\nu)=\sup_{\zeta\in C}\chi_{C-C}(\zeta+\nu)$. Then
$\beb\in\ac$, since $\beb$ is finitely supported.  Keeping
$\la\in\Lambda$ and $u\in C$,  we  estimate as follows:
\begin{eqnarray*}
\lefteqn{H(\la+u)=\sum _{\nu \in \Lambda } (\bbh\ast \alb  \ast
\alb ^* ) (\nu ) \,  \chi _{C-C}(\la+u-\nu) }\\&\le& \sum _{\nu
\in \Lambda } (\bbh\ast \alb  \ast \alb ^* ) (\nu ) \, \sup_{u\in
C} \chi _{C-C}(\la+u-\nu)\\&=&\sum _{\nu \in \Lambda } (\bbh\ast
\alb \ast \alb ^*) (\nu )\beb(\la-\nu)=(\bbh\ast \alb \ast \alb
^*\ast\beb) (\la ).
\end{eqnarray*}
Thus, $\sup_{u\in C}H(\la+u)\le (\bbh\ast \alb  \ast \alb
^*\ast\beb) (\la )$ and $\bbh\ast \alb  \ast \alb
^*\ast\beb\in \ac$. Consequently, 
$H\in W(\ac)$ and
\begin{equation}\label{eq:15}
\|H\|_{W(\ac)}\le \|\bbh\ast \alb  \ast \alb ^*\ast\beb\|_\ac\le
C\|\bbh\|_\ac,
\end{equation}
with $C= \|\alb  \ast \alb ^*\ast\beb\|_\ac$.
\end{proof}

Theorem~\ref{aldia} provides us with an almost diagonalization of
the operator $\sigma^w$ with symbol $\sigma\in \widetilde
M^{\infty,\Cal A}$ with respect to any  Gabor system $\Cal
G(g,\Lambda)$ that forms a tight frame.


The almost diagonalization result gives a characterization of the
symbol class $\widetilde M^{\infty,\ac}$ in two versions, a
continuous
one and a discrete one. We elaborate the discrete version. Let
$M(\sigma)$ be the matrix  with entries
\[
M(\sigma)_{\la,\mu}= \langle \sigma ^w \pi (\mu )g, \pi
  (\lambda )g\rangle \, , \quad \quad \lambda, \mu \in \Lambda \, .
\]
Clearly, $M(\sigma)$ depends on $g$, and we assume that the window
satisfies the assumptions of Theorem~\ref{aldia}.

\begin{theorem} \label{almost2}
$\sigma \in \widetilde M^{\infty ,\cA}$
\fif\ $M(\sigma ) \in \cC _{\cA }$. Moreover,
\begin{equation}\label{eq:23}
c\| \sigma \|_{\widetilde M^{\infty,\ac}} \le \|M(\sigma )
\|_{\Cal C_\ac}\le \| \sigma \|_{\widetilde M^{\infty,\ac}} \, .
\end{equation}
\end{theorem}

\begin{proof}
In view of Theorem~\ref{aldia}, specifically $(i) \,
\Leftrightarrow \, (iii)$, we need to prove only the norm
equivalence. From estimate~(\ref{eq:8a}) it follows immediately,
that
\[
\|M(\sigma ) \|_{\Cal C_\ac}\le \|\bbh\|_\ac=\|\gras\circ
j\|_{W(\ac)}=\|\sigma\|_{\widetilde M^{\infty,\ac}}.
\]

To deduce  the converse inequality,  we use (\ref{eq:8}) to
conclude that $\|\gras\circ j\|_{W(\ac)}\le \|H\|_{W(\ac)}$ for
every $H$ having the dominating property of (\ref{eq:6}). In
(\ref{eq:14}) we have described an explicit dominating function
$H$ derived from an arbitrary sequence $\bbh$ satisfying
(\ref{eq:6b}). Clearly, we can choose a sequence $\bbh$ such that
$\|\bbh\|_\ac=\|M(\sigma ) \|_{\Cal C_\ac}$. Since
$\|H\|_{W(\ac)}\le C\|\bbh\|_\ac$ by (\ref{eq:15}), we obtain that
\[
\|\sigma\|_{\widetilde M^{\infty,\ac}}=\|\gras\circ
j\|_{W(\ac)}\le \|H\|_{W(\ac)}\le C\|\bbh\|_\ac=C\|M(\sigma )
\|_{\Cal C_\ac}.
\]

\end{proof}

The matrix $M(\sigma)$ is the matrix of the operator $\sigma^w$
with respect to the tight Gabor frame $\Cal G(g,\Lambda)$. Indeed,
let $f = \sum _{\la \in \Lambda } \langle f, \pi (\la)g\rangle \pi
(\la ) g$ be the tight frame expansion of $f\in \lrd $. It  has
the coefficients
\begin{equation}\label{stftl}
V_g^\Lambda f(\la)=\langle f, \pi (\la)g\rangle \, ,
\end{equation}
where  $V_g^\La f$ is just the restriction of  the STFT of $f$
to  the lattice $\La$. Clearly, one has
\[
  \langle \sigma^w f, \pi (\lambda )g\rangle
  = \sum _{\mu \in \Lambda } \langle f, \pi (\mu
  )g\rangle \, \langle \sigma^w \pi (\mu )g,  \pi (\lambda ) g
  \rangle \, ,
\]
or in other words
\begin{equation}\label{eq:18}
 V_g^\La (\sigma ^wf) =M(\sigma) V_g^\La f\,.
\end{equation}
This commutation relation can be depicted by the  diagram
\begin{equation} \label{diagram1}
\begin{matrix}
&\lrd   & \stackrel{\sigma ^w}{\longrightarrow} & \lrd  & \cr
&\downarrow {V_g^\La} & &\downarrow V_g^\La & \cr &\ell
^2(\Lambda)\, &\stackrel{M(\sigma )}{\longrightarrow} &\ell
^2(\Lambda ) &
\end{matrix}
\end{equation}

For a symbol  $\sigma\in\symbo $, the operator $\sigma ^w$ is
always bounded on $\lrd $, because $\cA \subseteq \ell ^1(\Lambda
)$ by Theorem~\ref{alla}, and thus  $\widetilde
M^{\infty,\ac}\subset M^{\infty,1}\subset \Cal B(\lrd)$. By
Sj\"ostrand's fundamental boundedness result every $\sigma ^w$
with $\sigma \in M^{\infty,1}$ is bounded on $\lrd $. Thus the
diagram makes perfect sense for $\sigma \in \symbo $.
Clearly, if $M(\sigma )$ is bounded on $\ell ^2(\Lambda )$, then
$\sigma ^w$ is bounded on $\lrd $. In particular, if $M(\sigma )$
satisfies Schur's test, then $\sigma ^w$ is bounded on $\lrd $ (see
also \cite{BC94} and 
\cite{sjo07}). However, in this case,  it is not clear  how to
recognize or characterize  the symbol of such an operator,  quite in contrast to
Theorem~\ref{aldia} 


If $g \in M^{1,1}$ and $\cG (g, \Lambda )$ is a Gabor frame, then
the range of $V_g^\Lambda $ is always a proper closed subspace of
$\ell ^2(\Lambda )$~\cite{fg97jfa}. Therefore the matrix $M(\sigma
)$ is not uniquely determined by the diagram~\eqref{diagram1}. The
next lemma taken from ~\cite{gro06} explains the additional
properties of $M(\sigma ) $.

\begin{lemma} \label{matrixprop}
  If $\sigma ^w $ is bounded on $\lrd $, then $M(\sigma )
  $ is bounded on $\ell ^2( \Lambda )$ and maps $\mathrm{ran}\, V_g^\La$
  into $\mathrm{ran}\, V_g^\La$ with $\mathrm{ker}\, M(\sigma ) \supseteq
  (\mathrm{ran}\, V_g^\La )^\perp$.

Let $T$ be a matrix such that $ V_g^\La (\sigma ^wf) =T V_g^\La f
$
for all $f\in \lrd $. If $\mathrm{ker}\, T \supseteq
(\mathrm{ran}\,
V_g\Lambda ) ^{\perp }$, then $T=M(\sigma )$.
\end{lemma}

Note that $M(\sigma ) $ is never invertible, because its kernel
always
contains the nontrivial subspace $  (\mathrm{ran}\, V_g^\La
)^\perp$.

\section{Properties of operators with symbols in $\widetilde M^{\infty,\Cal A}$}

In this section, we study  several properties and applications of
the symbol class $\symbo $. The common theme is how  properties of
the algebra $\cA $ are inherited by properties of operators with
symbols in $\symbo $.  We investigate the algebra property, the
inverse-closedness, and the  boundedness of \psdo s on distribution spaces.

\subsection{The Algebra Property.} Let us first consider the algebra
property. In Theorem~\ref{aldia} we have used consistently and
crucially that $\cA $ is an algebra. This property is inherited by
$\symbo $.

\begin{theorem}\label{sjowei}
$\widetilde M^{\infty,\ac}$ is a Banach algebra with respect
to the twisted product $\sharp$ defined in~\eqref{eq:ll1}.
\end{theorem}

\begin{proof}
We assume that $\cG (g, \Lambda )$ is a tight Gabor frame with $V_gg
\in W(\cA )$. 
Let $\sigma , \tau \in \widetilde M^{\infty,\ac}$. By
using \eqref{diagram1} several times,  we get that
\begin{eqnarray*}
  M(\sigma \, \sharp \, \tau ) V_g^\La f &=& V_g^\La ( (\sigma \, \sharp \,
  \tau )^w  \, f) = V_g^\La (\sigma^w \, \tau ^w f) \\
&=& M(\sigma ) \big( V_g^\La (\tau ^w f )) = M(\sigma ) M(\tau )
V_g^\La f \, .
\end{eqnarray*}
Therefore the operators $M(\sigma\, \sharp\, \tau )$ and $M(\sigma
) M(\tau )$ coincide on $\mathrm{ran} \, V_g^\La$. Since
$\sigma^w$, $\tau^w$ and $(\sigma\,\sharp\,\tau)^w$ are bounded on
$\lrd$, we can use Lemma~\ref{matrixprop} to conclude that both
$M(\sigma\, \sharp\, \tau )$ and $M(\sigma ) M(\tau )$ are zero on
$(\mathrm{ran} \, V_g^\La)^\perp$. Thus, we get the following
matrix identity
\begin{equation}
  \label{eq:21}
M(\sigma \, \sharp \, \tau ) =   M(\sigma ) M(\tau )\,.
\end{equation}
Since by Theorem~\ref{almost2} both $M(\sigma )$ and $M(\tau )$ are in $\cC
_{\cA  }$, the algebra property of $\cC _{\cA }$ implies that
$M(\sigma \, \sharp \, \tau ) \in \cC _{\cA }$. By Theorem~\ref{almost2} once
again, we deduce that $\sigma \, \sharp \, \tau \in \symbo $. 
 The  norm estimate follows from 
 $$
\|\sigma \, \sharp \, \tau \| _{\widetilde M^{\infty,\ac}} \leq C
\|M(\sigma \, \sharp \, \tau )\|_{\cca} \leq  C \, \|M(\sigma
)\|_{\cca}\, \|M(\tau )\|_{\cca} \leq C \, \|\sigma \|_{\widetilde
M^{\infty,\ac}} \, \|\tau \|_{\widetilde M^{\infty,\ac}} \, , $$
and so $\symbo$ is a Banach algebra.
\end{proof}

\subsection{Boundedness} Let $\cY$ be a solid space of sequences on the lattice $\La\subset
\rdd $. If $\ac$ acts boundedly on $\cY$ under convolution, then
it is
natural to expect that the
class $\widetilde M^{\infty,\ac}$  acts boundedly on a
suitable
function space associated to $\cY$.
The appropriate function spaces are the (generalized) \modsp s. We
present the definition of \modsp s that is most suitable for our
purpose.

\begin{definition}\label{genmod}
  Let $\cY$ be a solid Banach space of sequences on $\Lambda $ such that
the finitely supported  sequences are dense in $\cY$. Let $\cG (g,
\Lambda )$ be a tight frame for $\lrd $ such that $V_gg \in
W(\cA)$. Let $\cL _0$ be the span of all finite linear
combinations $f = \sum _{\lambda \in \Lambda } c_\lambda \pi
(\lambda )g$.

Now we define a norm on $\cL _0$ by
\begin{equation}
  \label{eq:hh3}
  \|f\|_{M(\cY)}=\|V_g^\La f\|_\cY \, .
\end{equation}
The \modsp\ $M(\cY)$ is the norm completion of $\cL _0$ with
respect to
the $M(\cY)$-norm.
\end{definition}

Using this definition, the following general boundedness theorem
is very easy.

\begin{theorem}\label{bound}
Let $\cA $ be a solid involutive Banach algebra with respect to
convolution and let $\cY$ be a  solid  Banach space of sequences
on
$\Lambda $.   If  convolution of $\cA $  on $\cY$ is bounded, $\cA
\ast
\cY\subseteq \cY$,  and if  $\sigma \in
\moda$, then $\sigma ^w $ is bounded on $M(\cY)$. The operator
norm
can be estimated uniformly by $$ \|\sigma^w \|_{M(\cY) \to M(\cY)}
\leq
\|M(\sigma )\|_{\cca} \le \|\sigma \|_{\moda} \, .$$
\end{theorem}
\begin{proof} For $f\in\cL _0$ we use the commutative diagram  (\ref{eq:18}) and
(\ref{eq:ma3}) and obtain the boundedness in a straightforward
manner
from the estimate
\begin{eqnarray*}
\|\sigma^wf\|_{M(\cY)}&=&\|V_g^\La(\sigma^wf)\|_\cY=\|M(\sigma)V_g^\La
f\|_\cY\\
&\le & \|M(\sigma)\|_{\cca}\|V_g^\La f\|_\cY=
\|M(\sigma)\|_{\cca}\|f\|_{M(\cY)}.
\end{eqnarray*}
Since $M(\cY)$ is defined as the  closure of $\cL _0$, the above
estimate extends to
all $f\in M(\cY)$.
\end{proof}

 The proof looks rather trivial, but it hides many issues. In fact,
 our definition of the \modsp\ $M(\cY)$ is based on the main result
 of \tfa\  and coorbit theory. The general definition of a \modsp\
 ~\cite[Ch.~11.4]{book} starts with the solid function space $W(\cY )$
 on $\rdd $. The \modsp\ $M(\cY )$ then consists of all distributions
 $f$ such that $V_gf\in W(\cY )$. If $\cY$ is a weighted sequence
 space $\lpm $ or $\ell ^{p,q}_m $ for $1\leq p,q < \infty $, one
 obtains the standard \modsp s $\Mmpq $ discussed in Section~2. The main result of \tfa\ shows that the continuous definition can
 be discretized. If $\{ \pi (\lambda ) g: \lambda \in \Lambda \}$ is a
 frame for $\lrd $ with a ``nice'' window $g$, then $f\in M(\cY )$
 \fif\ $V_g^\Lambda \in \cY $; and this is the definition we are
 using. This characterization also implies the independence of the
 symbol class $\symbo$ of the Gabor frame $\cG (g, \Lambda )$ in 
Definition~\ref{def:sjocl}. 

However, the above  characterization of \modsp s by means of Gabor
 frames involves the major results of \tfa\ and is far from trivial. See~\cite{fg89jfa,fg92chui}
 and \cite[Chs.~11-13]{book} for a detailed account. 

\subsection{Spectral Invariance}

Next we discuss the invertibility of \psdo s and the symbol of the
inverse operator. We would like to characterize the inverse of
$\sigma
^w$ in terms of the  matrix $M(\sigma )$. This point is a bit more
sublte, because the invertibility of $\sigma^w$ on $\lrd$ does not
guarantee the invertibility of $M(\sigma)$ on $\ell^2(\La)$ (see
Lemma~\ref{matrixprop}). To go around this difficulty, we need the
notion of a  pseudo-inverse.

Recall that an operator $A: \ell ^2 \to \ell ^2$ is
pseudo-invertible, if there exists a closed subspace $\cR
\subseteq \ell ^2$, such that $A$ is invertible on $\cR$ and
$\mathrm{ker}\, A = \cR ^\perp $. The unique operator $A^\dagger $
that satisfies $A^\dagger A h = A A^\dagger h = h $ for $h\in \cR$
and $\mathrm{ker}\, A^\dagger = \cR^\perp $ is called the
(Moore-Penrose) pseudo-inverse of $A$. The following lemma is an
important consequence of Theorem~\ref{bask} and is  taken
from~\cite{gro06}.

\begin{lemma}[Pseudo-inverses]   \label{pseudo}
If $\ac $ is inverse-closed in $\cB (\ell ^2)$   and $A\in \cca$
has a (Moore-Penrose)
pseudo-inverse $A^\dagger$, then $A^\dagger \in \cca $.
\end{lemma}

\begin{proof}
By means of the Riesz functional calculus~\cite{rudin73}  the
pseudo-inverse can be written as
$$
A^{\dagger }  = \frac{1}{2\pi
i} \int _{ C } \frac{1}{z}\, (z\mathrm{I}-A)\inv \, dz \, ,
$$
where $C$ is a suitable path surrounding $\mathrm{Sp}_{\cB (\ell
^2)} (A) \setminus \{ 0\}$.  Theorem~\ref{bask}   implies that
$\cC _{\cA }$ is inverse-closed in $\cB (\ell ^2)$. Hence 
$(z\mathrm{I}-A)\inv\in\cca$, and  by \eqref{eq:ca11} this formula
makes
sense in  $\cca$.  Consequently,  $A^\dagger \in \cca $.
\end{proof}

\begin{theorem}\label{wiener} If $\ac$ is inverse-closed in $\cB (\ell ^2(\Lambda
))$, then the class of pseudodifferential operators with symbols in  $\moda$ is
inverse-closed in $\cB (\lrd )$. In the standard formulation, if
$\sigma \in \symbo $ and $\sigma ^w $ is invertible on $\lrd $, then
$(\sigma ^w)\inv = \tau ^w$ for some $\tau \in \symbo $. 
\end{theorem}

\begin{proof}
Assume that $\weyl$ is invertible on $\lrd$ for some
$\sigma\in\moda$. Let $\tau \in \cS ' (\rdd )$ be the unique
distribution such that $\tau ^w = (\weyl )\inv $. We need to show
that $\tau\in\moda$.

Since $\tau^w$ is bounded on $\lrd$,  Lemma~\ref{matrixprop}
implies  that the matrix $M(\tau ) $ is
 bounded on $\ell ^2(\La) $ and  maps
$\mathrm{ran} \, V_g^\La $ into $\mathrm{ran} \, V_g^\La $ with
$\mathrm{ker} \, M(\tau ) \supseteq (\mathrm{ran}\, V_g^\La )
^{\perp }$.

If $f\in\lrd$, then by (\ref{eq:18})  we have 
$$
  M(\tau ) M(\sigma ) V_g^\La f = M(\tau ) V_g^\La (\weyl f) =   V_g^\La (\tau
  ^w \weyl f) = V_g^\La f \, .
$$
This means that $  M(\tau ) M(\sigma ) = \mathrm{Id}$ on
$\mathrm{ran}\, V_g^\Lambda $ and that $M(\tau ) M(\sigma ) = 0$
on
$(\mathrm{ran}\, V_g^\Lambda )^\perp $.
Likewise, $M(\sigma ) M(\tau ) = \mathrm{Id}_{\mathrm{ran}\,
  V_g^\La}$ and $\mathrm{ker}\, M(\tau )=(\mathrm{ran} \, V_g^\La )^{\perp} $.
Thus, we  conclude that $M(\tau ) = M(\sigma ) ^\dagger $.

By Theorem~\ref{almost2}, the hypothesis  $\sigma \in \moda $
implies that $M(\sigma )$ belongs to the matrix algebra $\cca$.
Consequently, by Lemma~\ref{pseudo}, we also have that $M(\tau ) =
M(\sigma )^\dagger \in \cca$. Using Theorem~\ref{almost2} again,
we conclude that $\tau \in \moda$.
\end{proof}

The proof of Theorem~\ref{wiener} should be compared to the proofs of
analogous statements by Beals~\cite{beals77} and
Sj\"ostrand~\cite{Sjo95}. The key element in our proof is Theorem~\ref{bask}
about the inverse-closedness of the matrix algebra $\cC _{\cA}$. The
hard work in Theorem~\ref{wiener} is thus relegated to the theory of Banach
algebras. 

\begin{example}\label{example3} Since the algebras $\ell^q_{v}(\Lambda )$,
$q=1,\infty$ of Examples~\ref{example1} and
\ref{example2} are inverse-closed, the corresponding classes of
pseudodifferential operators are inverse-closed as well. Using
(\ref{lpm}) we obtain that the class $M^{\infty,q}_{1\otimes
v\circ j^{-1}}(\rdd)$ is inverse-closed for $q=1,\infty$ and
suitable weights $v$. In particular,
$M^{\infty,1}_{1\otimes\lan\cdot\ran^s}(\rdd)$ is inverse-closed
for $s\ge 0$ and
$M^{\infty,\infty}_{1\otimes\lan\cdot\ran^s}(\rdd)$ is
inverse-closed for $s>2d$.
\end{example}

\begin{corollary}
[Spectral Invariance on Modulation Spaces] Let $\ac$ be an
inverse-closed algebra that acts boundedly on a solid space $\cY$.
If $\sigma ^w$ is invertible  on $\lrd $, then $\sigma ^w$ is
invertible on $M(\cY)$. That is,
\begin{equation}
  \label{eq:hh1}
   \mathrm{Sp}_{\cB(M(\cY))} (\weyl) \subseteq  \mathrm{Sp}_{\cB (\lrd)}
 (\weyl) \qquad \text{ for all } \, \sigma \in \moda \, .
\end{equation}
\end{corollary}

\begin{proof}
By Theorem~\ref{wiener}, we have that $(\weyl )\inv = \tau ^w$ for
some $\tau \in \moda$. Thus, Theorem~\ref{bound} implies 
that $\tau ^w$ is bounded on $M(\cY)$. Since $\weyl \tau ^w = \tau
^w
\weyl = \mathrm{I}$ on $\lrd$, this factorization extends to
$M(\cY)$. Therefore, $\tau ^w = (\weyl )\inv $ on $M(\cY)$.
Applied to the
operator $(\sigma -\lambda \mathrm{I})^w$ for $\lambda \not \in
\textrm{Sp}_{\cB (L^2)}$, we find that also $\lambda \not \in
\textrm{Sp}_{\cB (M(\cY))}$, and the inclusion of the spectra is
proved.
\end{proof}


\section{H\"ormander's Class and Beals' Functional Calculus}

The H\"ormander class $S_{0,0}^0(\rdd)$ consists of smooth
functions all of whose  derivatives are bounded,
\[
S_{0,0}^0(\rdd)=\{f\in C(\rdd):|D^\alpha f|\le C_{\alpha}\}.
\]
Clearly,  $S_{0,0}^0(\rdd)=\bigcap_{n\ge0}C^n(\rdd)$,
where $C^n(\rdd)$ is the space of functions with $n$ bounded
derivatives. A characterization with \modsp s was mentioned  by
Toft~\cite{toft07} (Remark 3.1 without proof) based on his results for  embeddings
between \modsp
s and Besov spaces. 
\begin{lemma}\label{horclas}
\begin{equation}
\bigcap _{n\geq 0} C^n (\rd
)=\bigcap_{s\ge0}M^{\infty,\infty}_{1\otimes\lan\cdot\ran^s}(\rd)=\bigcap_{s\ge0}M^{\infty,1}_{1\otimes\lan\cdot\ran^s}(\rd).
\end{equation}
Hence $S_{0,0}^0(\rdd) =
\bigcap_{s\ge0}M^{\infty,\infty}_{1\otimes\lan\cdot\ran^s}(\rdd)=\bigcap_{s\ge0}M^{\infty,1}_{1\otimes\lan\cdot\ran^s}(\rdd).$
\end{lemma}

\begin{proof}  We give a proof based on the technique
  developed in \cite{book}. We  shall prove three inclusions.

i) $\bigcap _{n\geq 0} C^n (\rd )   \subset
\bigcap_{s\ge0}M^{\infty,\infty}_{1\otimes\lan\cdot\ran^s}$. By
formula (14.38) of \cite{book} every  $f\in C^n(\rd)$ fulfills 
\[
\sup_{x\in\rd}|V_gf(x,\xi)|\le C \|f\|_{C^n} |\xi^\beta|^{-1},
\]
for all multi-indices  $|\beta|\le n$. Therefore, if
$f\in\bigcap_{n\ge0}C^n(\rd)$, then for every $\beta$ there is a
constant $C_\beta$ such that $ \sup_{x\in\rd}|V_gf(x,\xi)|\le C_\beta
|\xi^\beta|^{-1}.$ 
Since $\langle \xi \rangle ^n \leq \sum _{|\beta   | \leq n}
c_\beta
|\xi ^\beta |$ for suitable coefficients $c_\beta \geq 0$
(depending
on $n$), we obtain that, for every $n\geq 0$,
\begin{equation}
  \label{eq:hh2}
  \sup _{x,\xi \in \rd } \langle \xi \rangle ^n \, |V_g f (x,\xi )|
  \leq C_n <\infty \, ,
\end{equation}
whence $f \in
\bigcap_{s\ge0}M^{\infty,\infty}_{1\otimes\lan\cdot\ran^s}$.


ii) 
$
\bigcap_{s\ge0}M^{\infty,\infty}_{1\otimes\lan\cdot\ran^s}\subset
\bigcap_{r\ge0}M^{\infty,1}_{1\otimes\lan\cdot\ran^r}$. If
$s>r+d$,  then
\begin{eqnarray*}
\|f\|_{M^{\infty,1}_{1\otimes\lan\cdot\ran^r}}
&=&\int_{\rd}\,\sup_{x\in\rd}|V_gf(x,\xi)|{\lan\xi\ran^r}\,d\xi
\\ & \le & \sup _{x,\xi \in \rd } |V_g f (x,\xi ) | \langle \xi
\rangle ^s \,  \int_{\rd}\lan \xi\ran^{r-s}d\xi = C_{r-s} \,
\|f\|_{M^{\infty,\infty}_{1\otimes\lan\cdot\ran^s}} <\infty.
\end{eqnarray*}
The  embedding $M^{\infty,\infty}_{1\otimes\lan\cdot\ran^s}
\subset M^{\infty,1}_{1\otimes\lan\cdot\ran^r}$ for $s>r+d$ yields
  immediately the embedding   $
  \bigcap_{s\ge0}M^{\infty,\infty}_{1\otimes\lan\cdot\ran^s}\subset 
\bigcap_{r\ge0}M^{\infty,1}_{1\otimes\lan\cdot\ran^r}$.

iii) $
\bigcap_{s\ge0}M^{\infty,1}_{1\otimes\lan\cdot\ran^s}\subset
\bigcap _{n\geq 0} C^n(\rd )$. Let $f\in
\bigcap_{s\ge0}M^{\infty,1}_{1\otimes\lan\cdot\ran^s}$. Since
$V_gf(x,\xi)=(fT_xg)^\wedge (\xi)$, we may rewrite the
$M^{\infty,1}_{1\otimes\lan\cdot\ran^s}$-norm as
\begin{equation}\label{mi1}
\|f\|_{M^{\infty,1}_{1\otimes\lan\cdot\ran^s}}=\int_{\rd}\,\sup_{x\in\rd}|(fT_xg)^\wedge
(\xi)|{\lan \xi\ran^r}\,d\xi <\infty  \, .
\end{equation}
By (\ref{mi1}) and
some basic properties of the Fourier transform,  we majorize the
derivatives of $fT_xg$ by
\begin{equation}\label{mi2}
\|D^\alpha(fT_xg)\|_\infty\le
\|\big(D^\alpha(fT_xg)\big)^\wedge\|_1\le
C_\alpha\int_{\rd}|(fT_xg)^\wedge(\xi)|\,|\xi^\alpha|\,d\xi
\end{equation}
\[
\le
C_\alpha\int_{\rd}|(fT_xg)^\wedge(\xi)|\,\lan\xi\ran^{|\alpha|}\,d\xi
\le
C_\alpha\|f\|_{M^{\infty,1}_{1\otimes\lan\cdot\ran^{|\alpha|}}}\,
,
\]
and this  estimate is uniform in $x\in\rd$. 
Next, let $g\in\Cal S(\rd)$ be a function that is equal to $1$ on
the unit ball $B$. By the Leibniz rule we have
\[
D^\alpha(fT_xg)=\sum_{\beta\le \alpha}{\alpha\choose
\beta}D^\alpha f D^{\alpha-\beta}(T_xg)=D^\alpha f \, 
T_xg+\sum_{\beta< \alpha}{\alpha\choose \beta}D^\alpha f
D^{\alpha-\beta}(T_xg)
\]
If $y\in x+B$, then clearly  the sum with $\beta<\alpha$ vanishes
and
\[
D^\alpha (fT_xg)(y)=D^\alpha f\, T_xg(y)=D^\alpha f\, T_x\chi_B(y).
\]
Thus, by (\ref{mi2}), we get that
\[
\| D^\alpha f \|_\infty = \sup _{x\in \rd } \|D^\alpha
f\, T_x\chi_B\|_\infty \le\|D^\alpha(f\, T_xg)\|_\infty\le
C_\alpha\|f\|_{M^{\infty,1}_{1\otimes\lan\cdot\ran^{|\alpha|}}}.
\]
Consequently, if $f\in
\bigcap_{s\ge0}M^{\infty,1}_{1\otimes\lan\cdot\ran^s}$, then $f\in
\bigcap _{n\geq 0} C^n(\rd )$.
\end{proof}

The above lemma and Theorem~\ref{aldia} yield the following
characterization of the H\"ormander class.

\begin{theorem}\label{hordia}
Take a window $g\in\Cal S(\rd)$ such that $\Cal G(g,\Lambda)$ is a
tight Gabor frame for $L^2(\rd)$. Then the following are
equivalent
\begin{itemize}
\item[(i)] $\sigma\in S_{0,0}^0(\rdd)$.
\item[(ii)] $\sigma\in \Cal S'(\rdd)$ and for every $s\ge 0$ there is a
constant $C_s$ such that
\begin{equation} \label{eq:6c}
|\lan \sigma^{w}\pi(z)g,\pi(w)g\ran|\le C_s \lan w-z\ran^{-s}
\quad \text{ for }\, w,z\in\rdd.
\end{equation}
\item[(iii)] $\sigma\in \Cal S'(\rdd)$ and for every $s\ge 0$ there is a
constant $C_s$ such that
\begin{equation} \label{eq:6d}
|\lan \sigma^{w}\pi(\mu)g,\pi(\lambda)g\ran|\le C_s
\<\lambda-\mu\>^{-s} \quad \text{ for }\, \lambda,\mu\in\Lambda.
\end{equation}
\end{itemize}
\end{theorem}

\begin{proof} We apply  Theorem~\ref{aldia} with the
algebra $\Cal A=\ell^\infty_{\<\cdot\>^s}(\La)$ for $s>2d$. Since
$g\in\Cal S(\rd)$ \fif\ $V_gg \in \cS (\rdd )$~\cite{folland89,book}, we have
$V_gg\in W(\ac)$, and the relevant hypothesis on $g$ is satisfied.
Moreover, $\widetilde
{M}^{\infty,\ac}=M^{\infty,\infty}_{1\otimes\lan\cdot\ran^s}(\rdd)$
by (\ref{lpm}).

By  Theorem~\ref{aldia},   $\sigma\in
M^{\infty,\infty}_{1\otimes\lan\cdot\ran^s}(\rdd)$ for a given $s>
2d$, if and only if $$M(\sigma )_{\lambda \mu } = \langle \sigma
^w \pi (\mu )g,\pi (\lambda )g\rangle = \cO (\langle \lambda - \mu
\rangle ^{-s}).$$ Using Lemma~\ref{horclas},  we find that $\sigma
\in
S_{0,0}^0(\rdd)=\bigcap_{s\ge0}M^{\infty,\infty}_{1\otimes\lan\cdot\ran^s}(\rdd)$,
\fif\ $M(\sigma ) $ decays rapidly, as stated. The argument for
the equivalence $(i) \,  \Leftrightarrow \,(ii) $ is analogous.
\end{proof}


The Banach algebra approach to the H\"ormander class yields a new
proof of Beal's Theorem on the functional calculus in
$S_{0,0}^0(\rdd)$~\cite{beals77} (see also~\cite{Ueb88}).

\begin{theorem} If $\sigma \in S_{0,0}^0(\rdd)$ and $\sigma ^w$ is
  invertible on $\lrd $, then there exists a $\tau \in S_{0,0}^0(\rdd)$,
  such that $(\sigma ^w)\inv = \tau ^w$.
\end{theorem}
\begin{proof}
This follows from Lemma~\ref{horclas} and Example~\ref{example3},
where we have shown that the classes
$M^{\infty,1}_{1\otimes\lan\cdot\ran^s}(\rdd)$ and
$M^{\infty,\infty}_{1\otimes\lan\cdot\ran^s}(\rdd)$ are
inverse-closed, for $s\ge 0$ and $s>2d$ respectively.

\end{proof}

\section{Time-Frequency Molecules}


The result on almost diagonalization  can be
extended to much  more general systems than Gabor frames. These
are
so-called \tf\ molecules. They were introduced in~\cite{gro04},
a different version of \tf\ molecules was investigated
independently
in~\cite{BCHL06}.

\begin{definition} Let $\cG (g, \Lambda ) $ be a tight Gabor frame
  with $V_gg \in W(\cA )$. 
  We say that $\{e_\mu\in\rd
:\mu\in\Lambda\}$ forms a {\it family of $\cA $-molecules} if
there exists a sequence $\bba\in \ac$ such that for every
$\la,\mu\in\La$ we have
\begin{equation}
  \label{eq:cc23}
  |\lan e_\mu,\pi(\la)g\ran|\le \bba(\la-\mu)\, ,
\end{equation}
where the window $g$ satisfies the hypotheses of
Theorem~\ref{aldia}.
\end{definition}

In this context a set of \tfs s 
$\{\pi(\mu)g:\mu\in\La\}$ plays the role of ``atoms''.  Moreover, if
$V_gg\in W(\ac)$, then 
$V^\La_gg\in \ac$ and  $\{\pi(\mu)g:\mu\in\La\}$ is also a set of $\cA
$-molecules. 

The above definition does not depend on the choice of the window
$g$. Assume that $h$ also generates a tight frame $\cG (h,\Lambda )$ and
that $V_hh \in W(\cA )$.
Then 
$V_gh\in
W(\ac)$ as in the proof of Lemma~\ref{lem:htht}.  Thus, $V^\La_gh\in\ac$ and therefore the sequence $\bbb
(\lambda ) = |\lan h,\pi(\la)g\ran|$ is in $\cA $. 
We insert the tight frame expansion of $\pi (\lambda )h= \sum
_{\mu \in \Lambda } \langle \pi (\lambda ) h, \pi (\mu )g\rangle \pi
(\mu )g$ with respect
to $\cG (g,\Lambda )$ into~\eqref{eq:cc23}    
and obtain  that 
\begin{eqnarray*}
|\lan e_\mu,\pi(\la)h\ran|&\leq & \sum_{\nu\in\La}|\lan
\pi (\lambda )h,\pi(\nu)g\ran|\,|\lan e_\mu,\pi(\nu)g\ran| \\
&\le &
\sum_{\nu\in\La} \bba(\nu-\la)\bbb(\nu-\mu)=\bba^*\ast \bbb
(\la-\mu).  
\end{eqnarray*}
Since $\bba^*\ast \bbb\in \cA $, \eqref{eq:cc23} is also satisfied for
$h$ in place of $g$. 


 Although \tf\ molecules are rather different from
classical molecules in analysis that are defined by support and
moment conditions~\cite{stein-weiss,stein93}, they observe a similar principle. The
following statement  justifies the name ``molecules''.

\begin{proposition}\label{molbound}
Assume that $V_gg\in W(\cA )$.    If the pseudodifferential operator
$A$ maps the tight Gabor frame 
$\{\pi
(\lambda )g : \lambda \in \Lambda \}$ to a set of $\cA
$-molecules,
then $A $ is bounded on $\lrd $. More generally, $A$ is bounded on
every \modsp\ $M(\cY)$, whenever  $\cA $ acts continuously on $\cY $
by convolution. 
\end{proposition}

\begin{proof}
Let $\sigma \in \cS '(\rdd )$ be the distributional  symbol corresponding to $A$. 
We check the matrix of $\sigma ^w$ with respect to the Gabor frame
$\cG (g, \Lambda )$ and obtain 
$$
|\langle \sigma ^w \pi (\mu ) g, \pi (\lambda ) g\rangle | = |\langle
e_\mu , \pi (\lambda )g\rangle | \leq \mathbf{a}(\lambda - \mu ),
\qquad \qquad \lambda , \mu \in \Lambda \, .
$$
Since $\mathbf{a}\in \cA $, the equivalence (iii) $\Leftrightarrow
$(i) of Theorem~\ref{aldia}  implies that the symbol $\sigma $ of $A$
belongs to the class $\symbo$. Now Theorem~\ref{bound} shows the
boundedness of $A$ on $\lrd $ and on $M(\cY )$. 
\end{proof}

With the notion of molecules, we may rephrase Theorem~\ref{aldia}
as
follows.

\begin{corollary}
  Let $\ac$ be a solid involutive Banach algebra with  respect to
  convolution   on a lattice $\Lambda \subseteq \rdd $ and $\{\pi
  (\lambda ) g: \lambda \in \Lambda \}$ be a tight Gabor frame with a window
  $g$ such that $V_gg\in W(\ac)$. Then $\sigma
  \in \widetilde M^{\infty, \cA }$, \fif\ $\sigma^w$ maps the \tf\ atoms
  $\pi (\lambda )g, \lambda \in \Lambda ,$ into $\cA $-molecules.
\end{corollary}

 As a further  consequence we show that  Weyl operators with symbols
 in $\widetilde M^{\infty , \cA }$  are almost diagonalized with
   respect to $\cA $-molecules.

\begin{corollary}\label{chromole}
Let $\{e_\lambda:\lambda\in\Lambda\}$ and
$\{f_\mu:\mu\in\Lambda\}$ be two families of $\ac$-molecules. If
$\sigma\in \widetilde M ^{\infty,\Cal A}$, then there exists a
sequence $\tilde \bbh\in\Cal A$ such that
\begin{equation}\label{mole}
|\lan \sigma^{w}f_\mu,e_\lambda \ran|\le \tilde \bbh(\lambda-\mu)
\quad \text{ for } \lambda,\mu\in\Lambda.
\end{equation}
\end{corollary}
\begin{proof}
The argument is similar to the 
final implication in the proof of Theorem~\ref{aldia}.  Each
molecule $f_\mu$ has the tight frame expansion
\begin{equation*}
f_\mu = \sum _{\nu \in \Lambda } \langle f_\mu,\pi (\nu ) g
\rangle \pi (\nu ) g \, .
\end{equation*}
with $|\langle f_\mu , \pi (\nu )g\rangle | \leq \mathbf{a}'(\nu - \mu
)$ for some $\mathbf{a}'\in \cA $. Likewise for the molecules
$e_\lambda $. 
Therefore, for every $\la,\mu\in\La$ we have, by (\ref{eq:6b}),
that
\begin{eqnarray*}
\lefteqn{  |\langle \sigma ^w f_\mu,e_\la\rangle | = \Big|\sum
_{\nu, \nu' \in \Lambda } \, \langle f_\mu , \pi (\nu) g \rangle
\, \overline{\langle e_\la, \pi (\nu' )g\rangle} \, \langle \sigma
^w \pi (\nu)g, \pi (\nu' )g\rangle \Big|} \\ &\leq & \sum _{\nu,
\nu '\in \Lambda } \bba'(\nu- \mu) \bba( \nu'-\la) h (\nu'-\nu ) =
(\bba'\ast \bba^* \ast \bbh   )(\lambda -\mu ) .
\end{eqnarray*}
Since $\bba'\ast \bba^* \ast \bbh\in\ac$, we can take $\tilde
\bbh=\bba'\ast \bba^* \ast \bbh$ as the dominating sequence in
(\ref{mole}).
\end{proof}

We shall exhibit a family of $\ac$-molecules with polynomial decay in
the \tf\ plane. For this we take $\cA $ to be  the algebra
$\ac=\ell^\infty_{\lan \cdot\ran^s}(\Lambda)$   and $s>2d$.

\begin{lemma}\label{mole8}
Let $\ac=\ell^\infty_{\lan \cdot\ran^s}(\Lambda)$ and $s>2d$ and 
fix a tight Gabor frame window $g$ such that $V_gg\in W(\ac)$.
Now choose a set
$\{z_\mu:\mu\in\La\}\subset \rdd$ satisfying  $|z_\mu-\mu|\le C'$
for every $\mu\in\La$ and  a family $\{\varphi_{\mu}:\mu\in\La\}$ of
functions such that $|V_g\varphi_\mu(z)|\le C\lan z\ran^{-s}$ uniformly
in 
$\mu\in\La$.   Then the  collection
\[
e_{\mu}=\pi(z_\mu)\varphi_{\mu}\,,\quad \mu\in\La
\]
forms a family of $\ac$-molecules. 
\end{lemma}
\begin{proof}
For $\la,\mu\in\La$, the assumptions $|V_g\varphi_\mu(z)|\le C\lan
z\ran^{-s}$ and $|z_\mu-\mu|\le C'$ imply that
\[
|\lan
e_\mu,\pi(\la)g\ran|=|\lan\pi(z_\mu)\varphi_{\mu},\pi(\la)g\ran|=
|V_g\varphi_\mu(\la-z_\mu)| \le C\lan \la-z_\mu\ran^{-s}\le
C''\lan \la-\mu\ran^{-s},
\]
thus the assertion follows.
\end{proof}

As another application we show how \psdo s are almost diagonal with
respect to  local Fourier bases. To be specific, we consider the local
sine bases of the form
\[
\psi_{k,l}(t)=\sqrt{\frac 2 \alpha}\, b(t-\alpha k)
\sin\Big(\frac{(2l+1)\pi}{2\alpha}(t-\alpha k)\Big) \, .
\]
The ``bell function'' $b$ can be constructed to be real-valued and in
$C^N(\bR )$ for given smoothness $N\in \bN $, such that
the system $\{\psi_{k,l}:k\in\bz,
l=0,1,2,\dots\}$ forms an orthonormal basis of
$L^2(\br)$~\cite{auscher94,hernandez-weiss}. The
parametrization of all bell function starts from   a
function $\zeta \in C^{N-1}(\bR)$ that  is real,   even, with support in
$[-\epsilon,\epsilon]$ for  $\epsilon \in (0,\alpha/2)$ and satisfies $\int_\br
\zeta(s)ds=\frac\pi2$. Now set $ \theta(t)=\int_{-\infty}^t \zeta (s)\,ds$, then
the bell function $b$ is given by \[
b(t)=\sin(\theta(t))\cos(\theta(t-\alpha))\, .
\]
See~\cite{auscher94,hernandez-weiss} for the details of this construction.

 From this definition it follows that $b$ is a real-valued,  compactly supported
function in $C^N(\bR)$ that satisfies
$$
|V_bb(x,\xi )| \leq C \langle \xi \rangle ^{-N} \chi _{[-K,K]}(x) \leq
C \langle (x,\xi )\rangle ^{-N} \,
$$
for some constants $C,K>0$ (see, e.g., \cite{GS00} for a derivation of
this estimate). In particular, $b$ satisfies the standard assumption
$V_bb \in W(\ac)$ with respect to the
algebra $\cA = \ell ^\infty _{\langle   .\rangle ^N}(\La)$ and the lattice $\La=\alpha\bz\times\frac1{2\alpha}\bz$.

By splitting the sine into exponential and some algebraic
manipulations,  we may write the basis functions as
\begin{equation}\label{kompost}
\psi_{k,l}=\frac{(-1)^{kl}}{2i} \sqrt{\frac 2 \alpha}  \bigg(
\pi\Big(\alpha k,\frac 
l{2\alpha}\Big) b_+ -\pi\Big(\alpha k,-\frac
l{2\alpha}\Big)b_- \bigg),
\end{equation}
where  $b_{\pm} (t)=\, e^{\pm \pi i t/(2\alpha)} b (t)$. 
Since $|V_bb (z) | = \cO (\langle z \rangle ^{-N})$, we also have
$|V_{b_{\pm } } b_{\pm }(z) | = \cO (\langle z \rangle ^{-N})$ for all
choices of signs. Thus $V_{b_{\pm } } b_{\pm } \in W(\ac)$. Since the lattice $\La$ admits a tight Gabor frame window $g\in \Cal S(\br)$, we get that $V_gg\in W(\ac)$ and, therefore, $V_{g } b_{\pm } \in W(\ac)$. This allows us to apply Lemma~\ref{mole8} to conclude that $\{\pi (\la ) b_{\pm }:\la\in\La\}$ is a family of $\ac$-molecules.
By using the decomposition \eqref{kompost} and Corollary~\ref{chromole}, we obtain
the following form of almost diagonalization with respect to a local
sine basis. 

\begin{theorem}
Assume that $b\in W(\ell ^\infty _{\langle   .\rangle ^N})$ and
$\sigma\in M^{\infty,\infty}_{1\otimes \langle 
\cdot\rangle^s}(\br^2)$ for $N>s>2$.  Then there is a constant $C_s$
such that
\[
|\lan \sigma^{w}\psi_{k',l'},\psi_{k,l}\ran|\le C_s \bigg(
\bigg\lan
\Big(\alpha(k-k'),\frac1{2\alpha}(l-l')\Big)\bigg\ran^{-s} +
\bigg\lan
\Big(\alpha(k-k'),\frac1{2\alpha}(l+l')\Big)\bigg\ran^{-s} \bigg)
\]
for $k,k'\in\bz$ and $l,l'\in\bn\cup\{0\}$.
\end{theorem}
This is reminiscent of a result by 
Rochberg and Tachizawa in \cite{tachizawa-rochberg98}.
The same type of almost diagonalization holds for the other types of
local Fourier bases (where the sine is replaced by a pattern of sines
and cosines). 
\section{Appendix}

In this appendix we sketch the  proofs of Theorems~\ref{alla}
and~\ref{bask}.
By a basis change for the lattice $\Lambda $, we may assume
without
loss of generality that $\Lambda = \zd $.

First we show that $\ell ^1(\Lambda )$ is the \emph{maximal solid
  involutive convolution algebra over $\Lambda $} (Theorem~\ref{alla}).
 We prove
Theorem~\ref{alla}
by a sequence of elementary lemmas. Their common assumption is
that
$\cA $ is a solid involutive  Banach algebra with respect to
convolution over $\zd $.



\begin{lemma}\label{al2}
 $\ac $ is continuously  embedded into $\el$.
\end{lemma}
\begin{proof} If  $\bba\in \ac$, then  $\bba^**\bba\in \ac$ and by solidity
  $(\bba^* \ast \bba) (0) \delta _0 \in \cA $. On the one hand, we have
\[
(\bba^**\bba )(0)=\sum_{k\in\bz}\bba^*(-k)\bba(k)=\sum_{k\in\bz}\ov{\bba(k)}\bba(k)=\|\bba\|_2^2,
\]
on the other hand $$ \|\bba\|_2^2 \|\delta _0\|_{\cA } =
\|\bba^*\ast \bba(0) \delta _0 \|_{\cA } \leq \|\bba^* \ast \bba
\|_{\cA } \leq C \|\bba\|_{\cA } ^2 \, , $$ and thus $\cA $ is
continuously embedded in $\ell ^2(\zd )$.
\end{proof}

Let $\Cal B=\{\bba\in\el : \Cal F \bba\in L^\infty (\bT ^d)\}$
with the norm $\|\bba\|_{\cB  } = \|\cF \bba\|_\infty $, where the
Fourier transform is given by
\[
\Cal F \bba(\xi)=\sum_{k\in\zd }\bba(k)e^{-2\pi i k\cdot \xi}.
\]
Then $\Cal B$ is a Banach algebra. 

\begin{lemma}\label {ab}
$\ac $ is continuously embedded in $\Cal B$ and $\|\bba\|_{\Cal
B}\le \|\bba\|_\ac$ for $\bba\in\ac$.
\end{lemma}
\begin{proof}
Take an arbitrary $\bba\in\ac$. Denote the $n$-fold convolution
$\bba*\bba*\dots*\bba$ by $\bba^{*n}$. From Lemma~\ref{al2}  it follows that
\[
\|(\Cal F \bba)^n\|_2=\|\bba^{*n}\|_2\le C\|\bba^{*n}\|_\ac\le
C\|\bba\|^n_\ac,
\]
for some constant $C$. Thus, $\|(\Cal F \bba)^n\|_2^{1/n}\le
C^{1/n}\|\bba\|_\ac$ and by taking the limit as $n\to \infty$
we obtain the desired inequality $\|\Cal F
\bba\|_\infty=\|\bba\|_{\cB } \le \|\bba\|_\ac$.
\end{proof}


\begin{lemma}\label{domi}
If $\bbb\in\cB $ and $|\bba|\le \bbb$, then $\bba\in \cB $ and $\|
\bba\|_\cB\le \|\bbb\|_\cB$.
\end{lemma}
\begin{proof}
 Let $C_\bbc$ be the
convolution operator given by $C_\bbc \bbh=\bbh*\bbc$ acting on
$\bbh\in\el$. Then $C_\bbc $ is bounded on $\ell ^2(\zd )$ \fif\
$\cF \bbc \in L^\infty(\td) $,  and the operator norm on $\ell
^2(\zd )$ is precisely $\|C_\bbc\|_\text{op}=\|\Cal F
\bbc\|_\infty=\|\bbc\|_\cB$. If $|\bba| \leq \bbb$, then  we have
\[
|C_\bba \bbh|\le |\bbh|*|\bba|\le |\bbh|*\bbb=C_\bbb|\bbh| \qquad
\text{ for all } \bbh\in \ell ^2(\zd ) \, ,
\]
and so
\[
\|C_\bba \bbh\|_2\le \|C_\bbb|\bbh|\|_2\le
\|C_\bbb\|_\text{op}\|\bbh\|_2.
\]
Therefore,  $\|C_\bba\|_\text{op}\le\|C_\bbb\|_\text{op}$, that
is, $\|\bba\|_\cB\le \|\bbb\|_\cB$.
\end{proof}

\begin{lemma}\label{oneinf}
Let $\bba$ be a sequence on $\zd$ such that $\Cal F|\bba|$ is well
defined. Then
\begin{equation}\label{normeq}
 \|\bba\|_1=\|\Cal F |\bba|\|_\infty.
\end{equation}
\end{lemma}
\begin{proof}
 Assume first  that $\bba\in \ell ^1(\zd)$. In this case $\Cal
F|\bba|$ is continuous on $\td$, and therefore
\[
\|\Cal F |\bba|\|_\infty\ge \Cal F|\bba|(0)=\|\bba\|_1.
\]
Since  we always have   $\|\Cal F |\bba|\|_\infty\le
\|\bba\|_1$, the equality  (\ref{normeq}) holds.

Next  assume that $\Cal F|\bba|\in L^\infty(\td)$.  We  show that
$\bba\in \ell ^1(\zd)$, then (\ref{normeq}) follows. Consider the
truncation of $|\bba|$, $\bba_N=\ch_{[-N,N]^d}|\bba|$, where
$\ch_{[-N,N]^d}$ is the characteristic function of the cube
$[-N,N]^d \subset \bz ^d$ for  $N\in \bN$. Since
$|\bba_N|\le|\bba|$, Lemma~\ref{domi} yields that $\|\Cal F
\bba_N\|_\infty\le\|\Cal F|\bba|\|_\infty$ for all $N\in\bn$. As
$\bba_N\in \ell ^1(\zd)$, we obtain that
\[
\sum_{n\in[-N,N|^d}|\bba(n)|=\|\bba_N\|_1=\|\Cal
F|\bba_N|\|_\infty=\|\Cal F \bba_N\|_\infty\le\|\Cal
F|\bba|\|_\infty\,
\]
 for all  $N\in \bN$.  This shows that $\bba\in
\ell ^1(\zd)$.
\end{proof}

We now prove Theorem~\ref{alla} showing the maximality of  $\ell
^1(\zd)$ in the class of  solid involutive  convolution Banach
algebras over $\zd$.
\begin{theorem}\label{al1}
Assume that  $\Cal A$ is  a solid involutive Banach algebra of
sequences over $\zd $.  Then $\Cal A\subseteq \ell ^1(\zd)$ and
$\|\bba\|_1\le\|\bba\|_{\Cal A}$ for $\bba\in\ac$.
\end{theorem}
\begin{proof}
Let $\bba\in \cA$. Since  $\Cal A$ is solid, both  $\bba$ and
$|\bba|$ have the same
 norm in $\cA$. Thus Lemmas~\ref{al2} and~\ref{ab} imply that
\[
\|\bba\|_{\Cal A}=\||\bba|\|_{\Cal A}\ge \||\bba|\|_{\Cal
B}=\|\Cal F|\bba|\|_\infty \, ,
\]
and by Lemma~\ref{oneinf},  $\|\Cal F|\bba|\|_\infty =
\|\bba\|_1$. This shows that $\|\bba\|_{\Cal A}\ge\|\bba\|_1$ and
thus $\cA $ is embedded in $\ell ^1(\zd )$.
\end{proof}

\vspace{ 4 mm}

We next turn to the proof of Theorem~\ref{bask}. Let us recall the
basic facts about the Gelfand theory of  commutative Banach
algebras.

 (a) The convolution operator
$C_{\mathbf{a}}$ defined by $C_{\bba } \bbc = \bba \ast \bbc $ for
$\bba \in \ell ^1(\zd)$ is invertible, \fif\ the Fourier series
$\widehat{\mathbf{a}}(\xi ) = \sum _{k\in \zd } \bba(k) e^{2\pi i
k\cdot \xi } $  does not vanish at any $\xi \in \bT ^d$.

(b)  The Gelfand transform of $\mathbf{a}\in \ell ^1(\zd )$
coincides with the Fourier series $\widehat{\mathbf{a}}(\xi )$
for $\xi \in \bT ^d$.

(c) An element $a$ in a commutative Banach algebra is invertible,
\fif\ its Gelfand transform does not vanish at any point.

(d)  By Theorem~\ref{alla}, we have that $\cA \subseteq \ell
^1(\zd )$, therefore $\widehat{\ell ^1} \simeq \bT ^d \subseteq
\widehat{\cA }$, and the Gelfand transform of $\bba \in \cA $
restricted to $\bT ^d$ is just the Fourier series $\widehat{\bba
}$ of $\bba $.

\begin{proof}[Proof  of Theorem~\ref{bask}- First Part]
The equivalence of (i) and (iii)  follows  from the Gelfand theory
for commutative Banach algebras:

\textbf{(iii) $\, \Rightarrow \, $ (i).}  If $ \widehat{\cA }
\simeq \bT ^d$, then  the convolution operator $C_{\bba } $ is
invertible \fif\ $\widehat{\bba}(\xi ) \neq 0 $ for all $\xi \in
\bT ^d$ by (a), \fif\ $\bba $ is invertible  in $\cA $ by (b) --
(d).

\textbf{(i) $\, \Rightarrow \, $ (iii).} If $ \widehat{\cA } \not
\simeq \bT ^d$, then the invertibility criteria for $C_{\bba } $
in $\cB (\ell ^2)$ and for $\bba $ in $\cA $ differ, and $\cA $
cannot be inverse-closed in $\cB (\ell ^2)$.


\textbf{(ii) $\, \Rightarrow \, $ (i).} The convolution operator
$C_{\bba }$ has the  matrix $A$ with entries $A_{kl}=\bba(k-l)$.
Consequently $\bbd_A (l) = \sup _{k} |A_{k,k-l}|= |\bba(l)|$ and
$\|A\|_{\cca } = \|\bbd_A\|_{\cA } = \|\bba \|_{\cA }$. Thus $A
\in \cca $ \fif\ $\bba \in \cA $.

If $\cA $ is not inverse-closed in $\cB (\ell ^2)$, then there
exists an $\bba \in \cA $, such that $C_{\bba } $ (with matrix
$B$)  is invertible on $\ell ^2 (\zd )$  with inverse $C_{\bbb }
$, but $\bbb \not \in \cA $. This means that $B$ cannot be in
$\cca $ and so $\cca $ is not inverse-closed in $\cB (\ell ^2)$.

\textbf{(i) $\, \Rightarrow \, $ (iv).} Since $\cA $ is inverse-closed
in $\cB (\ell ^2)$, the spectrum of the convolution operator
$C_{\mathbf{a}}$ acting on $\ell ^2$ coincides with the spectrum of
$\mathbf{a}$ in the algebra $\cA$, $\mathrm{Sp}_{\cA}(\mathbf{a}) =
\mathrm{Sp}_{\cB (\ell ^2)} (C_\mathbf{a})$. In particular, the
spectral radii of the particular elements $\delta _\lambda $ and of
$C_{\delta _\lambda}$ coincide. On the
one hand we have
$$
r_{\cA }(\delta _\lambda) = \lim _{n\to \infty } \|\delta _\lambda
^n\|_\cA ^{1/n} =  \lim _{n\to \infty } \|\delta _{n\lambda}
\|_\cA ^{1/n} = \lim _{n\to \infty } \|\omega (n\lambda ) 
^n\|_\cA ^{1/n}  \, .
$$
On the other hand, the convolution operator $C_{\delta _\lambda }$ on
$\ell ^2$ is simply  translation by $\lambda $ and is
unitary. Therefore the spectral radius of $C_{\delta _\lambda } $ is
$r_{\cB (\ell ^2)} (C_{\delta _\lambda  }) = 1$. By inverse-closedness
  we obtain the GRS-condition
   $\lim _{n\to \infty } \|\omega (n\lambda ) 
^n\|_\cA ^{1/n}   = 1$ for all $\lambda \in \Lambda $. 

\textbf{(iv) $\, \Rightarrow \, $ (iii).}  Since $\mathbf{a}\in \cA $
has the expansion $\mathbf{a} = \sum _\lambda \mathbf{a}(\lambda )
\delta _\lambda $, we have $\|\mathbf{a}\|_{\cA} \leq \sum _\lambda
|\mathbf{a}_\lambda \|\delta _\lambda \|_{\cA} = \|\mathbf{a}\|_{\ell
  ^1 _\omega} $. Thus $\ell ^1 _\omega  $ is continuously embedded in
  $\cA $, and by Theorem~\ref{alla} $\cA \subseteq \ell ^1$. By a theorem
  of Gelfand-Naimark-Raikov~\cite{gelfandraikov} the spectrum of $\ell
  ^1_\omega $ is isomorphic to $\bT ^d$  \fif\ $\omega $ satisfies the
  GRS-condition. If (iv) holds, then we obtain the embeddings $\bT ^d
  \simeq \ell ^1_\omega \subseteq \widehat{\cA } \subseteq \ell ^1
  \simeq \bT ^d$. Thus $\cA \simeq \bT ^d$ as claimed. 
\end{proof}

The substance of the theorem lies in the implication (iii) $\,
\Rightarrow \, $ (ii). To prove this non-trivial implication, we
study the Fourier series associated to an infinite matrix $A=
(a_{kl})_{k,l \in \zd }$~\cite{deL75,Bas90}.

Let $D_A(n)$ be the $n$-th diagonal of $A$, i.e., the matrix with
entries
\begin{eqnarray*}
D_A(n)_{kl} &= &
\begin{cases}
a_{k,k-n} \quad \text{ if } l=k-n \\ 0   \quad \text { else }
\end{cases}
\end{eqnarray*}
for $k,l,n \in \zd $. Furthermore, define the ``modulation'' $M_t,
t\in \bT ^d$,  acting on a sequence $\bbc = (\bbc(k))_{k\in \zd }$
by
\begin{equation}\label{eq:ma5}
\big( M(t) \bbc \big)(k) = e^{2\pi i k\cdot t } \bbc(k) \, .
\end{equation}
Finally to every matrix $A$ we associate the matrix-valued
function
\begin{equation}\label{eq:ma6}
\mathbf{f}(t) = M_t A M_{-t} \qquad  t\in \bT ^d \, .
\end{equation}

Clearly each $M(t) $  is unitary on $\ell ^2 (\zd )$ and $M_{t+k}
=M_t $ for all $k\in \zd $. Consequently $\mathbf{f}(t) $ is $\zd
$-periodic and $\mathbf{f}(0) = A$. Furthermore,  $A $ is invertible on
$\ell ^2 (\zd ) $ \fif\ $\mathbf{f}(t) $ is invertible for all
$t\in \bT
^d$.

It is natural to study the Fourier coefficients and the Fourier
series of the matrix-valued function $\mathbf{f}(t)$. It has the
following properties~\cite{Bas90}.

\begin{lemma}\label{lem:vvfourier}
\begin{itemize}
\item[(i)] $\mathbf{f}(t) _{kl} = a_{kl} e^{2\pi i (k-l)\cdot t}$ for $k,l \in \zd
$, $t\in \td $.
\item[(ii)] The matrix-valued Fourier coefficients of $\mathbf{f}(t)$ are given by
\begin{equation}\label{eq:ma7}
\widehat{\mathbf{f}}(n) = \int _{[0,1]^d} \mathbf{f}(t) e^{-2\pi i n \cdot
t} \, dt
= D_A(n)
\end{equation}
(with the appropriate interpretation of the integral).
\item[(iii)] Let $ \cA (\td , \cB (\ell ^2) )$ be the space of all matrix-valued
Fourier expansions $\mathbf{g}$ that are given by $\mathbf{g}(t) =
\sum _{n\in \zd } B_n \, e^{2\pi i n\cdot t}$ with $B_n \in \cB
(\ell ^2)$ and $\big( \|B_n \|_{op} \big)_{n\in \zd } \in \cA $.
Then
\begin{equation}\label{eq:ma8}
A \in \cca \quad \Longleftrightarrow \quad \mathbf{f}(t) \in  \cA
(\td ,
\cB (\ell ^2) ) \, .
\end{equation}
\end{itemize}
\end{lemma}

\begin{proof}
(i) follows from a simple calculation. For (ii) we interpret the
integral entrywise and find that
\begin{eqnarray*}
\widehat{\mathbf{f}}(n) _{kl} &=& \int _{[0,1]^d} \mathbf{f}(t) _{kl}\,
e^{-2\pi i
n\cdot t}\, dt \\ &=& a_{kl} \int _{[0,1]^d}  e^{-2\pi i (n+l-k)
\cdot t} \, dt = a_{k,l} \delta _{n+l-k} \, ,
\end{eqnarray*}
and so $\widehat{\mathbf{f}}(n) = D_A (n)$ is the $n$-th side diagonal, as
claimed.

(iii) follows from the definition of $\cca $. According to (ii),
the $n$-th Fourier coefficient of $\mathbf{f}(t)  = M_tA M_{-t}$
is just the $n$-th diagonal $D_A(n)$ of $A$ and $\|D_A(n)\|_{op}=
\bbd_A(n)$. So if $A \in \cca $, then $\mathbf{f}(t) = M_tA
M_{-t}$ is in $\cA (\td , \cB (\ell ^2) )$. Conversely, if
$\mathbf{f}(t)  = M_tA M_{-t}\in \cA (\td , \cB (\ell ^2) )$, then
$(\bbd_A(n))_{n\in \zd } = \big( \|\hat{A}(n)\|_{op}\big)_{n\in
\zd } \in \cA $ and thus $A \in \cca $.
\end{proof}

To prove the non-trivial implication of Theorem~\ref{bask}, we
need Wiener's Lemma for matrix-valued functions.

Baskakov's proof makes use of  the Bochner-Phillips version of
Wiener's Lemma for absolutely convergent Fourier series with
coefficient in a Banach algebra~\cite{BP42}. The proof exploits  the relation
between the ideal theory and the  representation theory of a
Banach algebra,  and the description of invertibility by means of
ideals.

Let $\cA $ be a Banach algebra with an identity and $\cM \subseteq
\cA $ be a closed left ideal, i.e., $\cA \cM \subseteq \cM $. Then
$\cA $ acts on the Banach space $\cA / \cM $ by the left regular
representation
\begin{equation}\label{eq:ma9}
\pi _{\cM } (\bba) \tilde{x} = \widetilde{\bba x} \qquad \text{
for } \bba\in \cA, \tilde{x} \in \cA /\cM\, ,
\end{equation}
where $\tilde{x} $ is the equivalence class of $x$ in $ \cA /
\cM $.

The following lemmata are standard and can be found in any
textbook on Banach algebras, see e.g.,
~\cite{Rick60,BD73}.

\begin{lemma}\label{lem:irred}
If $\cM  $ is a \emph{maximal} left ideal of $\cA $, then $\pi
_{\cM }$ is algebraically irreducible. This means that the
algebraically generated subspace $\{ \pi (\bba) \tilde{x}: \bba\in
\cA \}$ coincides with $ \cA/ \cM $ for \emph{all} $\tilde{x} \neq
0$.
\end{lemma}

\begin{lemma}\label{lem:invert}
Let $\cA$ be a Banach algebra with identity. An element $\bba\in
\cA $ is left-invertible (right-invertible), \fif\ $\pi _{\cM
}(\bba) $ is invertible for every maximal left (right) ideal  $\cM
\subseteq \cA $.
\end{lemma}

\begin{lemma}[Schur's Lemma for Banach Space Representations]\label{lem:schur}
Assume that $\pi : \cA \to \cB (X)$ is an  algebraically
irreducible representation of $\cA $ on a Banach space $X$. If
$T\in \cB (X)$ and $T\pi (\bba) = \pi (\bba) T$ for all $\bba\in
\cA $, then $T$ is a multiple of the identity operator
$\mathrm{I}_X$ on $X$.
\end{lemma}

\begin{proof}[Proof of Theorem~\ref{bask} - Second Part]
\textbf{(iii) $\, \Rightarrow \, $ (ii).} Assume that $A \in \cca
$ is invertible in $\cB (\ell ^2)$. Then the associated  function
$\mathbf{f}$ defined by $ \mathbf{f}(t) = M_t A M_{-t}$ possesses
a $\cB (\ell ^2)$-valued Fourier series
\begin{equation}\label{boldf}
\mathbf{f}(t) = \sum _{n\in \zd } D_A(n) e^{2\pi i n\cdot t} \, ,
\end{equation}
 where $D_A (n)$ is the $n$-th side diagonal of $A$ and the
sequence $\bbd_A (n) = \|D_A (n)\|_{op}$ is in $\cA $.

We identify the commutative algebra $\cA $ with a
subalgebra of $\cA (\td , \cB (\ell ^2) )$ via the embedding
$j:\cA \to \cA (\td , \cB (\ell ^2) )$

\begin{equation}\label{eq:9a}
j(\bba ) (t) = \sum _{n\in \zd } \bba(n) \, e^{2\pi i n\cdot t}
\mathrm{I} = \widehat{\bba } (t) \mathrm{I} \, .
\end{equation}
Clearly, since $j(\bba )$ is a multiple of the identity operator
$\mathrm{I}$, $j(\bba )$ commutes with every $T \in \cA (\td , \cB
(\ell ^2) )$.


Now let $\cM $ be a maximal left ideal of $ \cA (\td , \cB (\ell
^2))$ and $\pi _{\cM }$ be the corresponding representation. Since
every $j(\bba )$ commutes with every $T\in \cA (\td , \cB (\ell
^2))$,  we  find that $$\pi _{\cM } (T ) \pi _{\cM } ( j( \bba )
)= \pi _{\cM } ( j( \bba ) ) \pi _{\cM } (T ) \qquad \forall T\in
\cA (\td , \cB (\ell ^2)), \bba \in \cA \, . $$

As a consequence of Lemma~\ref{lem:schur} on the algebraic
irreducibility of $\pi _{\cM}$, $\pi _{\cM } (j(\bba ))$ must be a
multiple of the identity, and since $\pi _{\cM }$ is a
homomorphism, there exists a multiplicative linear functional
$\chi \in \widehat{\cA }$, such that $\pi _{\cM }(j(\bba)) = \chi
(\bba ) \, \mathrm{I}$.

Here comes in the crucial hypothesis that $\widehat{\cA } \simeq
\td $. This hypothesis says that  there exists a $t_0 \in \td $
such that $\chi (\bba ) = \widehat{\bba }(t_0)$. Consequently,

\[
\pi _{\cM } (j(\bba )) =  \widehat{\bba }(t_0)\, \mathrm{I} \qquad
 \text{for all     } \bba \in \cA \, .
\]

Let $\delta _n$ be the standard basis of $\ell ^1(\zd )$. Since
$\cA $ is solid, $\delta _n \in \cA $ and we have $j(\delta _n)(t)
= e^{2\pi i n\cdot t}\mathrm{I}$. Thus, the function $\mathbf{f}$ given in
(\ref{boldf}) can be
written as $\mathbf{f} = \sum _{n\in \zd } D_A(n) j(\delta _n)$.
Consequently we have
\begin{eqnarray}
\pi _{\cM } (\mathbf{f}) &=&   \pi _{\cM } \Big(\sum _{n\in \zd }
D_A(n) j(\delta_n) \Big) \notag \\ &=&  \sum _{n\in \zd } \pi
_{\cM } (D_A(n)) \pi _{\cM } (j(\delta _n))  \notag \\ &=&  \sum
_{n\in \zd } \pi _{\cM } (D_A(n))   e^{2\pi i n\cdot t_0}
\mathrm{I} \label{eq:ca14} \\ &=&  \pi _{\cM } \Big( \sum _{n\in
\zd } D_A(n) e^{2\pi i n\cdot t_0}\Big) \notag \\ &=& \pi _{\cM }
( \mathbf{f}(t_0)) \, , \notag
\end{eqnarray}

If $A\in \cca $ is invertible in $\cB (\ell ^2)$, then clearly
$\mathbf{f}(t) = M_t A M_{-t}$ is invertible in $\cB (\ell ^2)$
for every $t\in \td$ and consequently $\pi _{\cM }
(\mathbf{f}(t_0))$ is (left-) invertible for every maximal left
ideal in $\cA (\td, \cB (\ell ^2))$. By \eqref{eq:ca14} we find
that $\pi _{\cM } (f)$ is invertible for every maximal left ideal
in $\cA (\td, \cB (\ell ^2))$. Thus, by Lemma~\ref{lem:invert},
$\mathbf{f}(t)$ is invertible in the algebra  $\cA (\td, \cB (\ell
^2))$. This means that $\mathbf{f}(t)\inv $ possesses a Fourier
series $$ \mathbf{f}(t)\inv = M_t A\inv M_{-t} = \sum _{n\in \zd }
B_n \, e^{2\pi i n\cdot t} \, $$ with $(\|B_n\|_{op})_{n\in \zd }
\in \cA $. Since by Lemma~\ref{lem:vvfourier}(b) the coefficients
$B_n$ are exactly the side diagonals of $A\inv $,
Lemma~\ref{lem:vvfourier}(c) implies that  $A\inv \in \cca $. This
finishes the proof of Theorem~\ref{bask}.
\end{proof}

\begin{remark} Clearly Theorem~\ref{bask} also  works for arbitrary
  discrete abelian groups 
as index sets. 
\end{remark}

\def\cprime{$'$} \def\cprime{$'$}


\end{document}